%% file: main.tex
\def\mycmd{0} 
\if\mycmd0
\documentclass{article} 
\else
\documentclass{elsarticle} 
\fi
\input{preamble}

\begin{document}
\if\mycmd0
\input{arxiv-header}
\else
\input{els-header}
\fi

\section{Introduction}

    Cardiac mechanics is the main driver of the heart function, as it is in charge of pumping blood through the entire body, so its accurate description is fundamental to cardiac modeling. Deformation is initiated by the propagation of an electric potential, which induces the contraction of the cardiomyocytes, i.e. the cardiac muscle cells. In addition, the heart presents complex interactions during its motion, among which there is the friction with the pericardium, the orientation of the cardiac fibers that yield a twisting motion to the heart contraction, and the blood flow in its inner chambers, guided by the opening and closure of the valves, to name a few \cite{levick2013introduction}. 
    
    The heart can be modelled as a continuum undergoing large deformations, which results in a highly nonlinear system of partial differential equations (PDEs) formulated in a non-trivial geometry \cite{smith2004multiscale,Holzapfel2002489}. This problem is usually solved by means of a Newton method, which requires the iterated solution of the tangent problem, given by the Jacobian matrix \cite{wright1999numerical}. The resulting problem is large, with roughly a million degrees of freedom used for realistic simulations, which requires the use of iterative methods, most commonly a preconditioned GMRES method \cite{saad2003iterative}. Nonlinear elasticity is a very difficult problem in both theory and practice. There are no a-priori estimates (except in the case of small loads \cite{carstensen2004priori}) and therefore preconditioners need to be adjusted to each specific case, mostly starting from those that work best in linear elasticity (e.g. \cite{smith1992optimal, griebel2003algebraic}). In the current literature, the most successful preconditioners for this problem belong to two families: algebraic multigrid (AMG) \cite{xu2017algebraic} and domain decomposition \cite{toselli2004domain}.
    
    Multigrid preconditioners are well established and are generally used as black-box solvers, since they have a large amount of tuning parameters. Robust default settings are available in high-performance libraries such as HYPRE \cite{falgout2002hypre}, PETSc \cite{petsc-user-ref}, and Trilinos \cite{heroux2012new}. AMG preconditioners are widely used in the computational mechanics community, see for example \cite{el2010iterative,adams2002evaluation,brezina2006parallel,franceschini2019robust} and also \cite{AUGUSTIN2016622,Jiang2020,karabelas2022accurate} for cardiac specific studies. Domain decomposition preconditioners make use of a partition of the domain into subdomains, usually computed by an external library such as PT-SCOTCH \cite{chevalier2008pt} or METIS \cite{karypis1997parmetis}, and can be classified into overlapping and non-overlapping according to the number of layers of degrees of freedom shared by each processor \cite{toselli2004domain}. We focus on the Balancing Domain Decomposition by Constraints (BDDC) preconditioner \cite{dohrmann2003preconditioner}, spectrally equivalent (\cite{li2006feti}) to the Finite Elements Tear and Interconnect -- Dual Primal (FETI-DP) preconditioner \cite{farhat2001feti}. These methods are two-level, non-overlapping, substructuring preconditioners where continuity across subdomain boundaries is enforced only at a small set of primal degrees of freedom.  BDDC has already been shown to provide good results in cardiac mechanics for the case of structured meshes and lowest order finite elements \cite{pavarinoSZ2015,colli2018numerical}. We remark that these previous studies focused only on the standard two-level BDDC algorithm, without considering multilevel variants.
    
    The main novelty of this work is to provide a systematic numerical study of the performance of multilevel BDDC preconditioners through its main parameters, and to compare it with a state-of-the-art AMG solver on both structured and unstructured finite element meshes, with linear and quadratic finite elements. More specifically, we make use of the BDDC implementation available in PETSc \cite{zampini2016pcbddc}, that we have embedded into the Deal.II and LifeX libraries \cite{dealii,lifex}, and the HYPRE interface for the boomerAMG preconditioner \cite{falgout2002hypre}. The development of HYPRE has received much more attention in the community, with many collaborators throughout its long development and improvement. This imbalance with respect to BDDC, developed by Stefano Zampini \cite{zampini2016pcbddc}, is dealt with in this work by fine-tuning its parameters so that the comparison is as fair as possible.

    This work is structured as follows. In Section \ref{section:mechanics}, we present the cardiac mechanics model and the discretization choices that we use. In Section \ref{section:bddc}, we briefly introduce the AMG and BDDC preconditioners and their relevant components, most importantly the primal/dual degrees of freedom and the local/coarse problems for BDDC. Details on how to prepare the preconditioners for optimal performance are given in Section \ref{section:paralimpl}, with code snippets that show the main developmnents done for this work. The benchmark tests considered are reviewed in Section \ref{section:benchmarks}. We compare several modern direct solvers for the local and coarse problems in Section \ref{section:tuning}, namely PETSc LU, SuperLU \cite{li2005overview}, KLU \cite{davis2010algorithm}, UMFPACK \cite{davis2007umfpack}, Mumps \cite{amestoy2000mumps}, and Intel\textsuperscript{\textcopyright} MKL-Pardiso, as well as the performance of different configurations of primal spaces. In Section \ref{section:scalability}, we study the strong scalability of the preconditioner and compare it to the AMG scalability in two scenarios, one with 1 million degrees of freedom (DoFs) solved with up to 128 CPUs, and the other with 8 million DoFs solved with up to 1024 CPUs. We conclude our work with the assessment of the performance of the preconditioners in a realistic electromechanical simulation considering the electrophysiology and systemic circulation in Section \ref{section:realistic}.

\section{Cardiac mechanics}\label{section:mechanics}
Mathematical models used for cardiac mechanics are based on the  Continuum Mechanics framework \cite{Holzapfel2002489}. We consider a reference domain $\Omega$ representing a human left ventricle with its boundary $\partial \Omega$ divided into three non-overlapping regions: basal $\partial\Omega_\text{base}$, epicardium $\partial\Omega_\text{epi}$, and endocardium $\partial\Omega_\text{endo}$. From the displacement $\displacement:\Omega\to \R^3$ we define the velocity $\dot\displacement=
\frac{\mathrm{d}\,\displacement}{\mathrm{d}t}$, the strain tensor $\tenF = \ten I + \nabla \displacement$ with $\ten I$ the identity tensor, its determinant $J = \det \tenF$ and the Green-Lagrange tensor $\ten E=\frac 1 2\left(\tenF^T\tenF - \ten I\right)$. The conservation of linear momentum is given by
    \begin{equation}\label{eq:momentum conservation}
    \ddot{\displacement} - \dive \ten P(\tenF) = \vec 0,
    \end{equation}
where $\ten P$, the Piola stress tensor, is decomposed into passive and active contributions \cite{quarteroni2017integrated,ambrosi2011electromechanical} as
    $$ \ten P(\tenF) = \ten P_\text{pas}(\tenF) + \ten P_\text{act}(\tenF). $$
The passive part is obtained from a Helmholtz potential $\Psi$ as 
    $$ \ten P_\text{pas} = \parder{\Psi}{\tenF} $$
and represents the mechanical response of the tissue. Throughout the benchmark tests we consider the Guccione potential \cite{Guccione1991}, given by a quasi-incompressible fibre-oriented exponential material:
\begin{align*}
   \Psi(\tenF) &= \frac C 2 \left(\mathrm{e}^{Q(\tenF)} - 1\right) + \frac B 2 (J-1)\log J, \\
   Q(\tenF) &= b_{ff}E_{ff}^2 + b_{ss}E_{ss}^2 + b_{nn}E_{nn}^2 + b_{fs}(E_{fs}^2 + E_{sf}^2) + b_{fn}(E_{fn}^2 + E_{nf}^2) + b_{sn}(E_{sn}^2 + E_{ns}^2), 
\end{align*}
where $E_{ab} = \ten E\vec a\cdot \vec b$ for $a,b\in\{f,s,n\}$ denote the components of $\ten E$ in the fiber-induced frame of reference $\{\vec f, \vec s, \vec n\}$. See \cite{Guccione1991} for reference values of the related parameters.

The active part instead represents the contraction of the cardiomyocytes induced by the propagation of an electric potential \cite{franzone2014mathematical}, and we consider it as
    $$ \ten P_\text{act} = \gamma(t)\frac{(\tenF \vec f)\otimes \vec f}{|\tenF \vec f|}, $$
where $\gamma$ is a given function known as the activation function. More details on the generation of the fibers and the activation function can be found in \cite{Bayer20122243} and \cite{franzone2014mathematical} respectively. 
\paragraph{Boundary conditions} On the epicardium $\partial\Omega_\text{epi}$ we consider a generalized Robin condition which takes into account the effect of the pericardium on the cardiac wall \cite{Usyk2002}: 
\begin{subequations}\label{eq:bc}
\begin{equation}\label{eq:bc epi}
\ten P(\tenF)\vec N = -\vec g_\text{epi}(\displacement, \dot \displacement) \coloneqq - (\vec N\otimes \vec N)\left(K_\perp^\text{epi}\displacement + C_\perp^\text{epi}\dot\displacement\right) - (\ten I-\vec N \otimes \vec N)\left(K_\|^\text{epi}\displacement + C_\|^\text{epi}\dot\displacement\right)=\vec 0\quad\text{ on }\partial\Omega_\text{epi},
\end{equation}
where $\vec N$ denotes the outward normal vector in reference configuration. On the endocardium $\partial\Omega_\text{endo}$ we consider the action of blood in the ventricle chamber, which reads
\begin{equation}\label{eq:bc endo}
   \ten P(\tenF)\vec N = -\vec g_\text{endo}(\tenF) \coloneqq  -p J\tenF^{-T}\vec N \quad\text{ on }\partial\Omega_\text{endo}. 
\end{equation}
Finally, on the ventricle base we consider a null traction: 
\begin{equation}
    \ten P(\tenF)\vec N = -g_\text{base}(\tenF) \coloneqq \vec 0 \quad\text{ on }\partial\Omega_\text{base}.
\end{equation}
\end{subequations}

\paragraph{The mechanics problem}
Putting together equation \eqref{eq:momentum conservation}, the boundary conditions \eqref{eq:bc} and the initial conditions $\displacement(0) = \displacement_0, \dot\displacement(0)=\vec v_0$ we obtain the weak formulation of the cardiac mechanics problem: Find a displacement $\displacement$ such that 
    \begin{equation}\label{eq:mechanics weak formulation}
        \int_\Omega \ddot{\displacement}\cdot \displacement^*\,dx + \int_\Omega \ten P(\tenF):\nabla \displacement^*\,dx
        + \int_{\partial\Omega}\vec g(\displacement, \dot \displacement, \tenF)\cdot \displacement^*dS = 0
    \end{equation}
for all test functions $\displacement^*$, where we used the definition
    $$ \vec g(\displacement, \dot \displacement, \tenF) \coloneqq \vec g_\text{epi}(\displacement, \dot \displacement)\bbone_\text{epi} + \vec g_\text{endo}(\tenF)\bbone_\text{endo} + \vec g_\text{base}(\tenF)\bbone_\text{base} $$
with $\bbone_\omega$ given by the index function, which takes the values 1 in $\omega$ and 0 in its complement. 

\subsection*{Numerical approximation}
{\bf Time discretization.} We consider a backward Euler time discretization with a constant time-step $\Delta t$, such that $t_i = t_0 + i\Delta t$ with the standard notation $\eta^n\approx \eta(t_n)$, for any time-dependent function $\eta$. \\
{\bf Space discretization.} We consider a triangulation of $\Omega$ into mesh-regular hexahedra \cite{QuarteroniValli2008} together with $\vec H^1$ conforming finite elements of orders 1 and 2, which we will denote with $\Q_1$ and $\Q_2$ respectively. The resulting nonlinear problem can be stated as: At instant $t_n$, given two previous displacements $\displacement^{n-1}$ and $\displacement^{n-2}$, find a displacement $\displacement^n$ such that 

\begin{multline}\label{eq:mechanics time discrete}
    \mathcal{F}(\displacement^n)\coloneqq  \int_\Omega \left( \frac{\displacement^n - 2 \displacement^{n-1} + \displacement^{n-2}}{\Delta t} \right)\cdot \displacement^*\,dx + \int_\Omega \ten P(\tenF^n):\nabla \displacement^*\,dx
        + \int_{\partial\Omega}\vec g\left(\displacement^n,  \frac{\displacement^n - \displacement^{n-1}}{\Delta t}, \tenF^n\right)\cdot \displacement^*dS = 0,
\end{multline}
where $\ten F(\displacement^n)\coloneqq \ten I + \nabla \displacement^n$.  The nonlinear problem \eqref{eq:mechanics time discrete} is solved by a Newton method: starting from an initial guess $\displacement^{n,0}\coloneqq \displacement^{n-1}$, the next Newton iteration is computed by solving the linear system 
\begin{equation}\label{eq:tangent problem}
    D\mathcal F(\displacement^{n,k-1})[\delta \displacement^{n,k}] = -\mathcal F(\displacement^{n,k-1}),
\end{equation}
for the displacement increment $\delta\displacement^{n,k}$,
and then by updating $\displacement^{n,k} = \displacement^{n,k-1} + \delta\displacement^{n,k}$. Here, $D\mathcal F(\vec a)[\vec b]$ stands for the Frechét derivative of $\mathcal F$ evaluated in $\vec a$ in the direction $\vec b$ \cite{ambrosetti1995primer}. At each Newton step,
the linear system (\ref{eq:tangent problem}) with coefficient matrix $D\mathcal F(\vec d^{n,k-1})$  is solved by a GMRES method with the BDDC/AMG preconditioner detailed in Section \ref{section:bddc}.

\section{The preconditioners}\label{section:bddc}
In this section we review the preconditioners under consideration, i.e. BDDC and AMG. A more detailed explanation of BDDC will be provided in order to give the relevant concepts to be tested in the parameter tuning section (Section \ref{section:tuning}).

\subsection{AMG preconditioner}
For a detailed review of algebraic multigrid methods, we refer to \cite{stuben2001review}. The main idea is that of solving the problem only approximately on a fine mesh, and then reducing the problem to another one on a coarser mesh. The approximate solution has the objective of reducing the finer components of the error, and then the coarse components are reucrsively interpreted as fine components on the coarse mesh. The algebraic element resides on the choice of the mesh, which is inferred from the matrix structure. Consider a linear system $\mat A\vec x = \vec b$, we can then define the strongly connected nodes in the matrix, or coarse nodes, as 

    $$ C = \left\{ J: j\in J \text{ if } A_{ij}\geq \varepsilon A_{ii} \right\},$$
for a given constant $\varepsilon$. Defining an interpolation operator $\mat P^T_C$ that goes from the coarse nodes $C$ to the original ones, we can compute the coarse operator

    $$ A_C = \mat P_C \mat A \mat P_C^T, $$
on which the same procedure can be repeated. This procedure is repeated until the coarse matrix $\mat A_C$ is small enough for it to be inverted easily with a direct method. After the projection phase, an extension phase if also performed, where the current iterate is interpolated using subsequent operators $\mat P_C^T$ until the original mesh resolution is obtained. The default configuration of HYPRE is extremely robust for scalar problems, but for elasticity there are two considerations to make: the first one is that the block-structure of the problem--given by the components of the displacement-- needs to be taken into consideration. The second one is that low energy modes, in this case rigid body motions, must be preserved by the projection operator. More details on this in Section \ref{section:paralimpl}.

\subsection{BDDC preconditioner}

This is a non-overlapping domain decomposition preconditioner, and can be seen as an evolution of the balancing Neumann-Neumann methods. Throughout this section, we adopt the notation from \cite{toselli2004domain} and follow \cite{colli2018numerical}, so we consider a decomposition of $\Omega$ into $N$ non-overlapping subdomains $\Omega_i$ of diameter $H_i$, $H\coloneqq\max_i H_i$, such that $ \Omega = \cup_i \Omega_i$ together the interface among them $\Gamma\coloneqq \cup_i\partial\Omega_i\setminus \partial\Omega$. By reordering the degrees of freedom according to the interior and interface with subscripts $I$ and $\Gamma$ respectively, the inverse of a matrix $\mat A$ can be factorized as 

$$ \mat A^{-1} = \begin{bmatrix}\mat I & -\mat A^{-1}_{II}\mat A_{I\Gamma}\\ \mat 0 & \mat I\end{bmatrix}        \begin{bmatrix}\mat A_{II}^{-1} & \mat 0 \\ \mat 0 & \mat S_{\Gamma}^{-1}\end{bmatrix}                 \begin{bmatrix}\mat I & \mat 0 \\ -\mat A_{\Gamma I}\mat A_{II}^{-1} \mat I\end{bmatrix}, $$
where we consider the Schur complement $\mat S_\Gamma\coloneqq \mat A_{\Gamma \Gamma} - \mat A_{\Gamma \mat I} \mat A_{II}^{-1} \mat A_{\Gamma I}^{-1}$ and $\mat A_{UV}$ denotes a sub-matrix according to the index sets $U$ (rows) and $V$ (columns). This factorization is exact, so the aim of the BDDC preconditioner is to provide a suitable approximation of $\mat S_\Gamma^{-1}$.

Now we define the relevant spaces used to build the preconditioner. For this, we denote with $V$ the discrete finite element space used to approximate the displacement, and by $V^{(i)}$ the local discrete space defined on $\Omega_i$\footnote{The space $V^{(i)}$ considers the local Dirichlet boundary conditions, but our formulation only has natural (Neumann and Robin) boundary conditions, so no modification is required. }. The local space is further split into the direct sum of its interior and interface subspaces as $V^{(i)}= V_I^{(i)}\bigoplus V_\Gamma^{(i)}$, which yields the global spaces 
    $$ V_I\coloneqq \prod_{i=1}^N V_I^{(i)}, \quad\text{ and }\quad V_\Gamma \coloneqq \prod_{i=1}^N V_\Gamma^{(i)}. $$
Note that functions in $V_\Gamma$ have possibly more than one value for each interface degree of freedom in $V$, meaning that functions in $V_\Gamma$ are not necessarily continuous across the subdomain interfaces. Indeed, this is a substructuring method 
\cite{toselli2004domain}, where the problem matrices remain unassembled throughout the entire solution process, meaning that there is no communication between processes to obtain the global matrix. We thus define the subspace
    $$ \widehat V_\Gamma \coloneqq \left\{ \text{functions of $V_\Gamma$ that are continuous across $\Gamma$}\right\},$$
as well as the intermediate space 
    $$ \widetilde V_\Gamma \coloneqq V_\Delta \bigoplus \widetilde V_\Pi, $$
where $\Pi$ and $\Delta$ stand for primal and dual degrees of freedom, $\widetilde V_\Pi$ corresponds to the subspace of functions which are continuous in a given set of degrees of freedom associated to the corners, edges and/or faces, and $V_\Delta\coloneqq \prod_{i=1}^N V_\Delta^{(i)}$ is the product space of the local subspaces $V_\Delta^{(i)}$ which vanish at the primal degrees of freedom. The choice of the primal degrees of freedom is fundamental for the construction of efficient and scalable BDDC preconditioners. It is well-known that considering only the vertices of the subdomains can give sub-optimal scalability, \cite{toselli2004domain}, but this problem can be alleviated by including the subdomain edge and/or face averages in the primal space.

We will further require restriction and extension operators represented by matrices with values in $\{0,1\}$:
    \begin{align*}
        \mat R_{\Gamma \Delta}: \widetilde V_\Gamma \to V_\Delta, \ \ \ \ \ \ \
        \mat R_{\Gamma \Pi}: \widetilde V_\Gamma \to \widehat V_\Pi, \\
        \mat R_\Delta^{(i)}: V_\Delta \to V_\Delta^{(i)}, 
        \ \ \ \ \ \ \
        \mat R_\Pi^{(i)}: \widehat V_\Pi\to \widehat V_\Pi^{(i)},
    \end{align*}
where $\widehat V_\Pi^{(i)}$ is the local subspace of primal interface functions. The passage from the unstructured spaces to $V$ (where the actual solution belongs to), requires the use of adequate averaging techniques in the extension operators. For the sake of exposition, we consider only the pseudo-inverse $\delta_i^\dagger(x)$ of the counting function defined at each degree of freedom $x$ on the interface of $\Omega_i$ by   
    $$ \delta_i^\dagger \coloneqq \frac 1{\text{ \#  \{subdomains sharing $x$\}}},$$
but note that more robust options could be the deluxe scaling [61] or anysotropic scalings as well [63,18,19]. Thus we define the local restriction operators $\mat R_{D, \Delta}^{(i)}$ by multiplying the nonzero element of $\mat R_\Delta^{(i)}$ by $\delta_i^\dagger$, and $\mat R_{D, \Gamma} \coloneqq \mat R_{\Gamma \Pi}\oplus \mat R^{(i)}_{D,\Delta}\mat R_{\Gamma\Delta} $. Finally, we denote with $\mat A^{(i)}$ the unassembled matrix corresponding to subdomain $\Omega_i$, where of course $\mat A = \sum_{i=1}^N\mat A^{(i)}$, and note that it can be written as 
    $$ \mat A^{(i)} = \begin{bmatrix} \mat A_{II}^{(i)} & \mat A_{I\Delta}^{(i)} & \mat A_{I \Pi}^{(i)} \\
        \mat A_{\Delta I}^{(i)} & \mat A_{\Delta\Delta}^{(i)} & \mat A_{\Delta\Pi}^{(i)}\\
        \mat A_{\Pi I}^{(i)} & \mat A_{\Pi \Delta} & \mat A_{\Pi \Pi}^{(i)} \end{bmatrix}. $$
With the previous definitions, the BDDC preconditioner (of $\mat S_\Gamma$) can be written as

    $$ \mat M_\text{BDDC}^{-1} = \mat R_{D,\Gamma}^T \widetilde{\mat S}_\Gamma^{-1}\mat R_{D,\Gamma},$$
where the approximate Schur complement $\widetilde{\mat S}_\Gamma$ is given by
    $$ \widetilde{\mat S}_\Gamma^{-1} = \mat R_{\Gamma \Delta}^T\left( \sum_{i=1}^N\begin{bmatrix}\mat 0 & \left(\mat R_\Delta^{(i)}\right)^T\end{bmatrix}\begin{bmatrix}\mat A_{II}^{(i)} & \mat A_{I\Delta}^{(i)} \\ \mat A_{\Delta I}^{(i)} & \mat A_{\Delta \Delta}^{(i)}\end{bmatrix}^{-1}\begin{bmatrix}\mat 0 \\ \mat R_\Delta^{(i)}\end{bmatrix} \right)\mat R_{\Gamma\Delta} + \Phi S_{\Pi\Pi}^{-1} \Phi^T,$$
and in addition we have
    \begin{align*}
        \mat S_{\Pi\Pi} &= \sum_{i=1}^N \left(\mat R_\Pi^{(i)}\right)^T\left(\mat A_{\Pi\Pi}^{(i)} - \begin{bmatrix}\mat A_{\Pi I}^{(i)} & \mat A_{\Pi \Delta}^{(i)}\end{bmatrix}\begin{bmatrix}\mat A_{II}^{(i)} & \mat A_{I\Delta}^{(i)} \\ \mat A_{\Delta I}^{(i)} & \mat A_{\Delta\Delta}^{(i)}\end{bmatrix}^{-1}\begin{bmatrix}\mat A_{I \Pi}^{(i)} \\ \mat A_{\Delta \Pi}^{(i)} \end{bmatrix}  \right)R_\Pi^{(i)}, \\
        \Phi &= \mat R_{\Gamma \Pi}^T - \mat R_{\Gamma \Delta}^T\sum_{i=1}^N \begin{bmatrix}\mat 0 & \left(\mat R_\Delta^{(i)}\right)^T\end{bmatrix}\begin{bmatrix}\mat A_{II}^{(i)} & \mat A_{I\Delta}^{(i)} \\ \mat A_{\Delta I}^{(i)} & \mat A_{\Delta \Delta}^{(i)}\end{bmatrix}^{-1}\begin{bmatrix}\mat A_{I\Pi}^{(i)} \\ \mat A_{\Delta \Pi}^{(i)}\end{bmatrix}\mat R_{\Pi}^{(i)}.
    \end{align*}
The columns of $\Phi$ represent the coarse basis functions, which are given by the minimum energy extension (with respect to the original bilinear form) of the primal constraints into each subdomain \cite{dohrmann2003preconditioner}.

As the number of subdomains increases, the coarse problem given by $\mat S_{\Pi \Pi}$ becomes the bottleneck. A recent solution for this is the approximation of such problem by means of a BDDC preconditioner, known as the \emph{multilevel} BDDC \cite{mandel2008multispace}, with the implementation details given in \cite{zampini2016pcbddc}. We note that this kind of solver is still under active research, where, to avoid an excessive deterioration in the problems's conditioning, adaptive coarse spaces can be computed by means of auxiliary eigenvalue problems \cite{mandel2012adaptive,pechstein2017unified}. The general outcome of this research, in albeit much simpler problems, is that multi-level strategies are fundamental to obtain optimal strong and weak scalability of the BDDC preconditioner in large scale problems.

\section{Parallel implementation}\label{section:paralimpl}
In order to achieve an adequate performance, the implementation of the BDDC interface in the deal.II library was done with special attention to memory efficiency. More specifically, exact memory preallocation was implemented in deal.II, together with a dedicated interface to the PETSc class \icode{PCBDDC}. We note that exposing all options in the preconditioner in not practical, as it is still possible to give a file with additional PETSc options through the \icode{-options\_file} command. The following routines have thus been included:

\paragraph{Matrix allocation.} Substructuring preconditioners such as BDDC require unassembled matrices, meaning that at no point there is a matrix assembly where shared degrees of freedom are added and scattered to their corresponding process. This is implemented in PETSc through a special type of matrix known as \icode{MATIS}, where each process owns a dedicated sequential matrix (\icode{MATSEQAIJ}). The locally owned degrees of freedom, extracted in deal.II with \icode{DofHandler::locally\_owned\_dofs}, and the degrees of freedom related to the locally owned elements, extracted with \icode{DoFTools::extract\_locally\_active\_dofs}, are the building blocks to compute then an exact local CSR structure given by the \icode{DynamicSparsityPattern} of the problem. The use of the latter ensures the scalability of the construction of the non-zero matrix entries. The resulting interface is easy to use, and allows for the reuse of the extracted degrees of freedom, as shown in Figure \ref{fig:dealii-interface-is}. Internally, \icode{reinit\_IS} leverages the problem connectivity to then use \icode{MatSeqAIJSetPreallocationCSR} from PETSc for the memory preallocation of the local matrix.

\begin{figure}[ht!]
\begin{lstlisting}[style=mystyle, language=C++]
// Sparsity pattern object, requires knowledge of local and
// ghosted degrees of freedom.
DynamicSparsityPattern dsp(locally_relevant_dofs);
// Matrix allocation uses the sparsity pattern and the 
// local degrees of freedom for both rows and columns.
system_matrix.reinit(locally_owned_dofs, 
                     locally_owned_dofs,  
                     dsp, 
                     mpi_communicator);
\end{lstlisting}
\begin{lstlisting}[style=mystyle, language=C++]
// For the unassembled matrix, additional information is required
// regarding the degrees of freedom from the local elements.
IndexSet locally_active_dofs;
DoFTools::extract_locally_active_dofs(dof_handler, locally_active_dofs);
// As before, the information is given separately for
// rows and columns. 
system_matrix.reinit_IS(locally_owned_dofs, locally_active_dofs,
                        locally_owned_dofs,ocally_active_dofs,
                        dsp,
                        mpi_communicator);
\end{lstlisting}
\caption{Comparison of the allocation of regular (top) and unassembled (bottom) matrices in the new interface for \icode{MATIS} matrices in deal.II.}
    \label{fig:dealii-interface-is}
\end{figure}

\paragraph{PCBDDC interface.} The interface for the preconditioner exhibits basic structure for setting the primal space and the symmetry of the problem. We have additionally included an interface for setting the coordinates of the degrees of freedom, which allows the preconditioner to obtain a better definition of the vertices, which can be ill-defined when handling unstructured grids. A simple usage of the preconditioner is shown in Figure \ref{fig:dealii-interface-bddc}.

\begin{figure}[ht!]
\begin{lstlisting}[style=mystyle, language=C++]
SolverControl solver_control(dof_handler.n_dofs(), 1e-12); // Stopping criteria
// Create CG solver and preconditioner as with any other.
LA::SolverCG solver(solver_control, mpi_communicator);
PETScWrappers::PreconditionBDDC preconditioner;
PETScWrappers::PreconditionBDDC::AdditionalData data;

// ------ Optional -------
// Here we compute the dof coordinates to improve the performance
// of the preconditioner. For vector problems (i.e. elasticity), it would also
// help setting the block size through the PETSc interface by
// using MatSetBlockSize(system_matrix.petsc_matrix(), dim).
std::map<types::global_dof_index, Point<dim>> dof_2_point;
DoFTools::map_dofs_to_support_points(MappingQ1<dim>(),
                                     dof_handler,
                                     dof_2_point);
IndexSet owned_dofs = dof_handler.locally_owned_dofs();
std::vector<double> coords(dim * owned_dofs.n_elements());
unsigned int i = 0; 
unsigned int d = 0;
for (auto it = owned_dofs.begin(); it != owned_dofs.end(); ++it)
{
    for (d = 0; d < dim; ++d)
        coords[dim * i + d] = dof_2_point[*it](d);
    ++i;
}
data.coords_cdim = dim;
data.coords_n = i;
data.coords_data = coords.data();
// ------ Optional end -------

// Initialize preconditioner and solve problem.
preconditioner.initialize(system_matrix, data);
solver.solve(system_matrix,
             solution,
             system_rhs,
             preconditioner);
\end{lstlisting}
\caption{Setup and usage of the BDDC preconditioner in deal.II.}
\label{fig:dealii-interface-bddc}
\end{figure}

\paragraph{Further improving the robustness of the preconditioners}
Both BDDC and AMG preconditioners are capable of using information about the structure of the problem. More specifically, the indices of the displacement are strided, meaning that they follow the pattern \texttt{\{x1, y1, z1, x2, y2, z2, \ldots, xN, yN, zN\}}. For BDDC, this is important as each component of a vertex can be considered as an independent primal degree of freedom, whereas for AMG this can be used to devise coarse nodes that follow this structure. This is defined with the \icode{MatSetBlockSize} PETSc command. 

In addition, low energy modes can be included through the \icode{MatSetNearNullSpace} function. For BDDC, this is important as rigid motions can be considered as primal degrees of freedom as well, whereas AMG uses them to restrain the projection operators from adding error frequencies associated with these modes. Note that to use this in BDDC, the \texttt{-pc\_bddc\_use\_nnsp} flag must be turned on.

\section{Benchmark tests}\label{section:benchmarks}
We perform all tests in three standard benchmarks \cite{wrro93701}: i) a loaded beam, ii) a swelling ventricle and iii) a contracting ventricle, which we describe in what follows. We highlight that the swelling and contraction tests are performed on unstructured hexahedral meshes. In all tests, both Q1 and Q2 finite elements have been employed. 

\begin{figure}[ht!]
        \centering
        \includegraphics[height=0.3\textwidth]{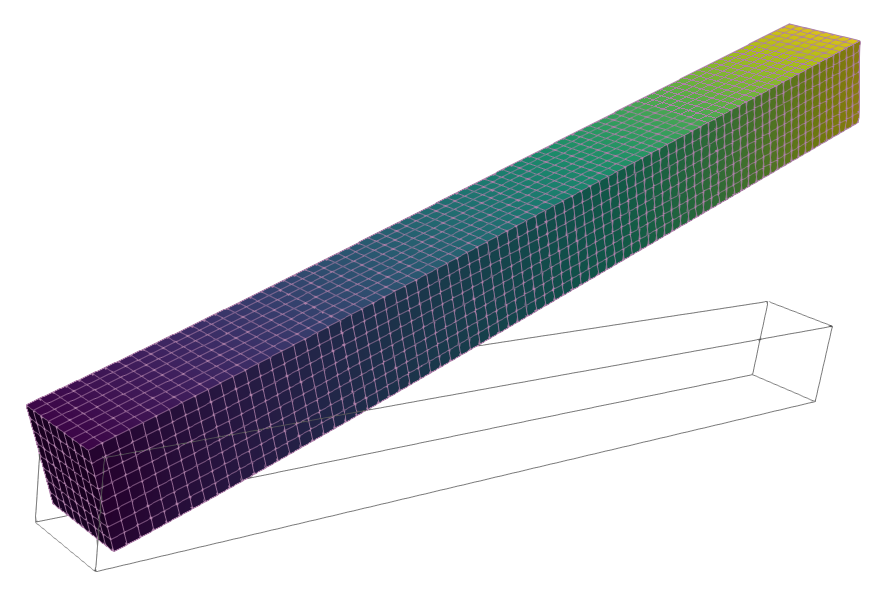}
        \caption{Beam solution at $t=0.1$.}
        \label{fig:beam}
    \end{figure}

\begin{figure}[ht!]
        \centering
        \includegraphics[height=0.4\textwidth]{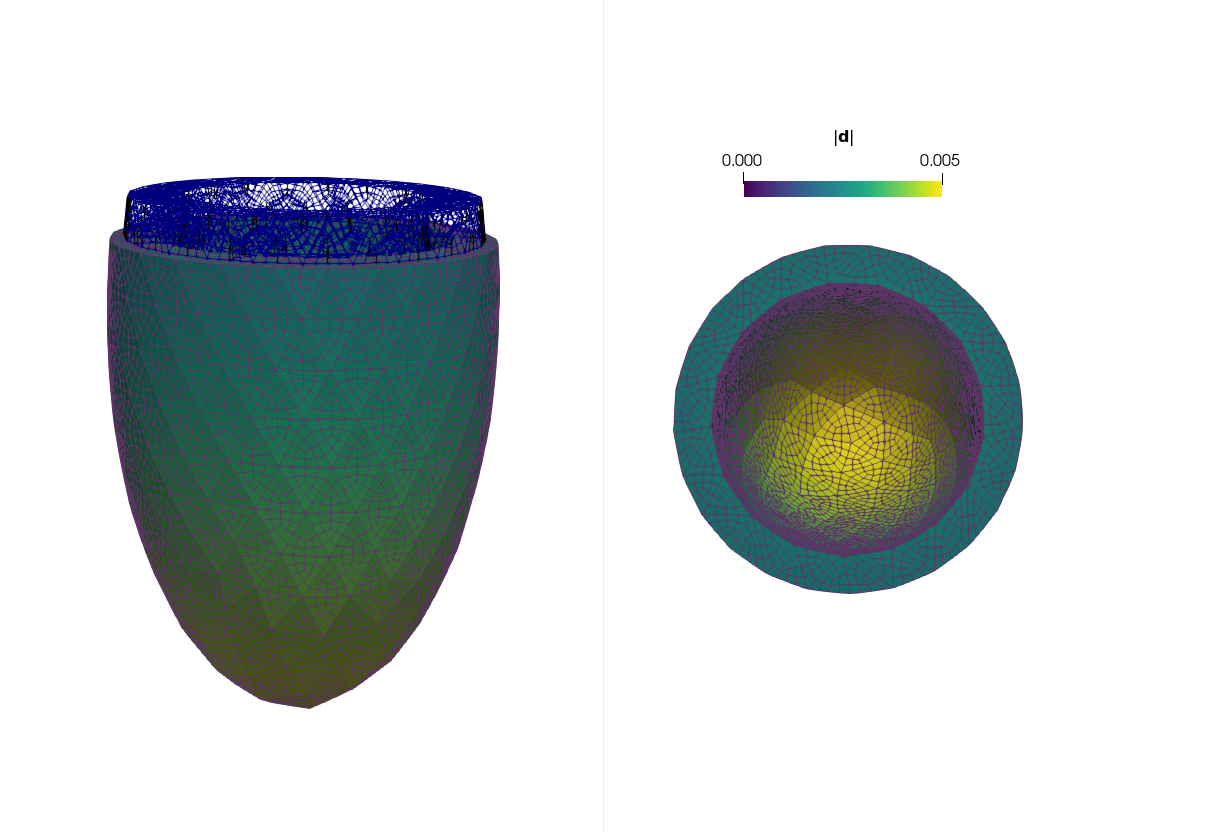}
        \caption{Swelling solution at $t=0.1$. Reference configuration is depicted in blue on the left.}
        \label{fig:swelling}
    \end{figure}

\paragraph{Beam test.} This test considers a structured mesh discretizing the domain $\Omega=(0,10)\times(0,1)\times (0,1)$ with $\partial\Omega_\text{endo}=\{z=0\}, \partial\Omega_\text{epi} = \{z=1\}$ and $\partial\Omega_\text{base}=\{x=0\}$, the rest undergoes null traction. Here we do not consider the active component, so $\ten P_\text{act}={\bf 0}$, and on the endocardium we consider an analytic pressure ramp given by $p(t) = 4\,\text{Pa}\, {\min\{t, 0.1\}}/{0.1}$, which at $t=1$ yield the deformation shown in Figure~\ref{fig:beam}. We use a time step of $\Delta t=10^{-3}$ seconds for this test.
\paragraph{Swelling test.} This test considers an unstructured grid, discretizing an idealized left ventricle, commonly referred to as a prolate ellipsoid. On the base we impose Robin boundary conditions with parameters $K_\perp^\text{base}=2\cdot10^5\,\text{Pa}/\text{m}, K_{\|}^\text{base}=2\cdot 10^4\,\text{Pa}/\text{m}, C_\perp^\text{base}=1 \cdot 10^4\,\text{s}\,\text{Pa}/\text{m}, C_{\|}^\text{base}=2\cdot 10^3\,\text{s}\,\text{Pa}/\text{m}$, whereas on the epicardium we impose a null traction condition \cite{pfaller2019}. We consider again a pressure ramp on the endocardium, given here by $p(t) = 2.5\,\text{kPa} {\min\{t, 0.1\}}/{0.1}$, which at $t=0.25$ yields the deformation shown in Figure~\ref{fig:swelling}. We use a time step of $\Delta t=10^{-3}$ seconds for this test.

\begin{figure}[ht!]
    \newcommand\wid{3cm}
        \centering
        \begin{subfigure}[b]{0.3\textwidth}
        \centering
        \includegraphics[height=\wid]{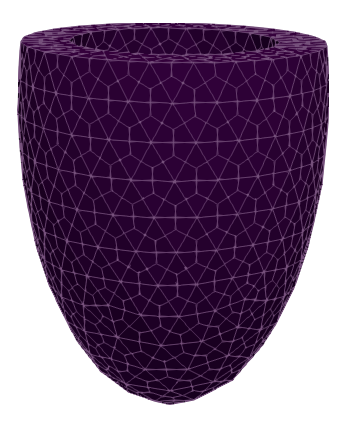}
        \caption*{$t=0.0$}
        \end{subfigure}
        \begin{subfigure}[b]{0.3\textwidth}
        \centering
        \includegraphics[height=\wid]{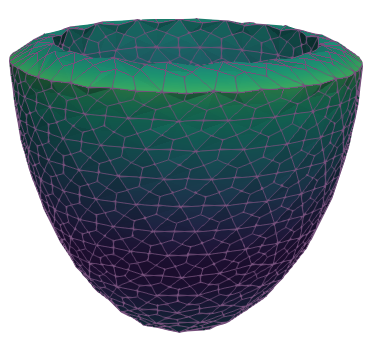}
        \caption*{$t=0.1$}
        \end{subfigure}
        \begin{subfigure}[b]{0.3\textwidth}
        \centering
        \includegraphics[height=\wid]{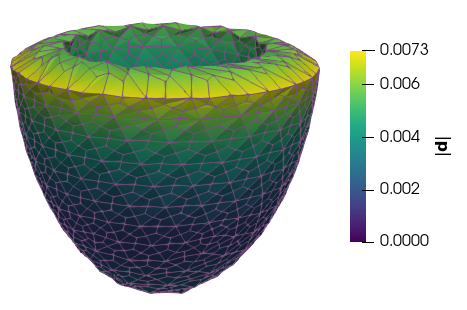}
        \caption*{$t=0.3$}
        \end{subfigure}
        
        \begin{subfigure}[b]{0.3\textwidth}
        \centering
        \includegraphics[height=\wid]{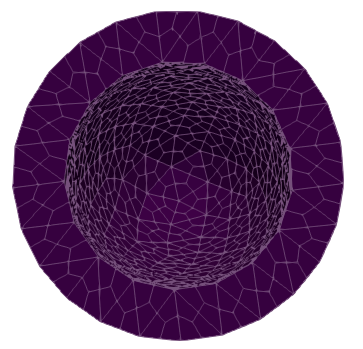}
        \caption*{$t=0.0$}
        \end{subfigure}
        \begin{subfigure}[b]{0.3\textwidth}
        \centering
        \includegraphics[height=\wid]{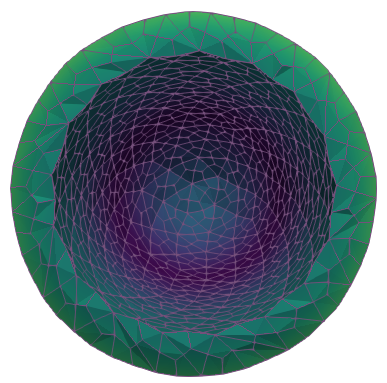}
        \caption*{$t=0.1$}
        \end{subfigure}
        \begin{subfigure}[b]{0.3\textwidth}
        \centering
        \includegraphics[height=\wid]{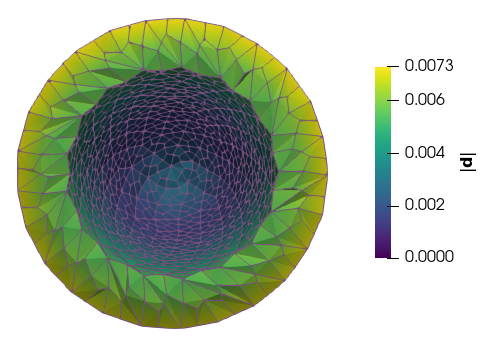}
        \caption*{$t=0.3$}
        \end{subfigure}
        \caption{Evolution of realistic case from front (top row) and above (bottom row).}
        \label{fig:contraction}
    \end{figure}

\paragraph{Contraction test.} In this test we consider the same geometry as in the swelling tests, but we include the Robin boundary conditions on both the base and the epicardium. We also consider the active stress term, with an analytic activation function given by $\gamma(t)=6\cdot 10^4 \min\{t, 0.25\}/0.25$ and an endocardial pressure given by $p(t)=15\,\text{kPa}\min\{t, 0.25\}/0.25$, which yields the deformation shown in Figure~\ref{fig:contraction}. We use a time step of $\Delta t=10^{-3}$ seconds for this test.

\section{BDDC tuning}\label{section:tuning}
In this section, we test two fundamental choices required by the BDDC preconditioner. The first one is the choice of direct solver for the local and coarse problems, and the second one is the choice of primal degrees of freedom for the coarse space.
\subsection{Choice of BDDC local and coarse direct solvers}
The choice of the different libraries used in this section is dictated by the ones available through PETSc. This results in testing PETSc LU, KLU, UMFPACK, SuperLU, Mumps and MKL--Pardiso for the serial direct solvers of the local problems and the libraries Mumps, SuperLU and MKL--Pardiso for the parallel direct solver of the coarse problem.

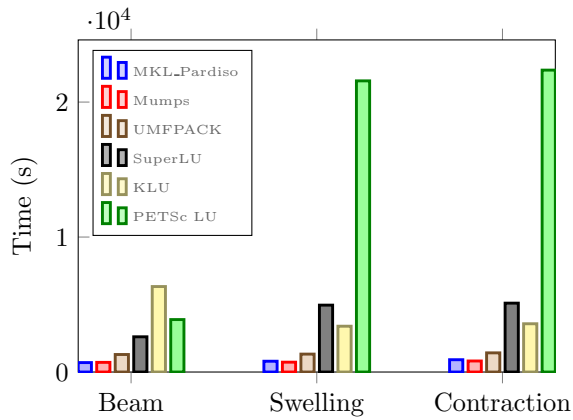
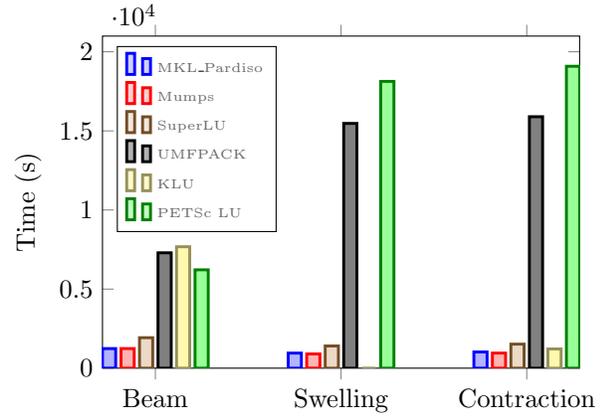
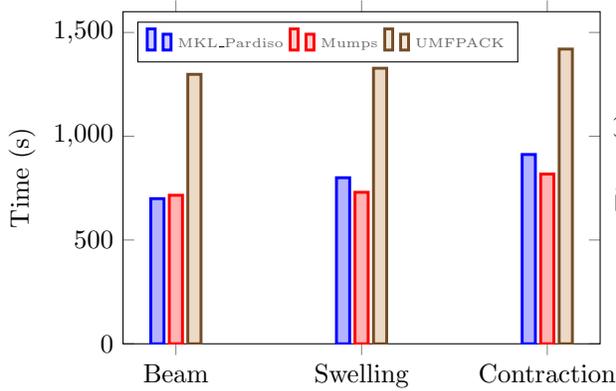
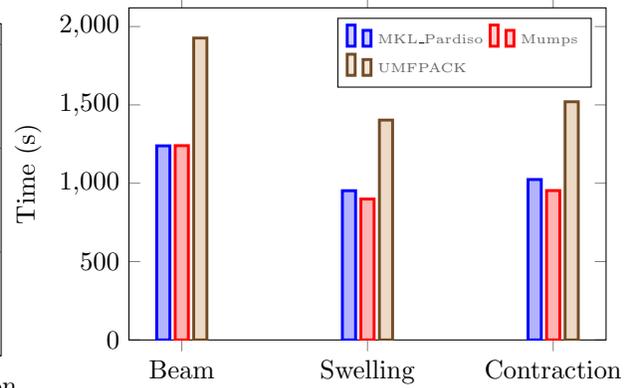
\begin{figure}
    \begin{subfigure}{0.48\textwidth}
        \centering
        \input{tikz/direct/local_solver_parallel_Q1}
        \caption{\Qone}
    \end{subfigure}
    \begin{subfigure}{0.48\textwidth}
        \centering
        \input{tikz/direct/local_solver_parallel_Q2}
        \caption{\Qtwo}
    \end{subfigure}
    
    \begin{subfigure}{0.48\textwidth}
        \centering
        \input{tikz/direct/local_solver_parallel_Q1_opt}
        \caption{\Qone \,(zoom on best 3 solvers)}
    \end{subfigure}
    \begin{subfigure}{0.48\textwidth}
        \centering
        \input{tikz/direct/local_solver_parallel_Q2_opt}
        \caption{\Qtwo \,(zoom on best 3 solvers)}
    \end{subfigure}
    \caption{Comparison of CPU time required with respect to the different local solver for first and second order finite elements. We compare the libraries MKL--Pardiso, Mumps, SuperLU, UMFPACK, KLU and PETSc LU. The second row shows the same results of the first one butzooming on the best 3 solvers Mumps, MKL--Pardiso and UMFPACK. }
    \label{fig:tuning:local}
\end{figure}

\paragraph{Local problem solvers.} We fix the Mumps library for the coarse problem. The model is solved with about 1 million degrees of freedom on 32 CPUs, which results in 15625 degrees of freedom per processor on average. The simulation is run for the first ten time steps. Note that the local problems are solved independently in each processor, meaning that this is indeed a serial test in each of the 32 processors. We show the results in Figure \ref{fig:tuning:local}, where we highlight that the 'swelling' and 'contraction' tests use unstructured meshes. In the first row we show the results for all libraries, where most notably PETSc LU, KLU and SuperLU have by far the worst performance. Indeed, for a clearer comparison we have removed them in the second row, where we show only MKL--Pardiso, Mumps and UMFPACK. In all cases, Mumps is the fastest direct solver by a small margin, except for the first order case in the 'beam' test, where the solution times of Mumps and MKL--Pardiso are comparable. Besides the increase in CPU time between first and second order finite elements, there is no qualitative difference in the results obtained from different finite element orders.

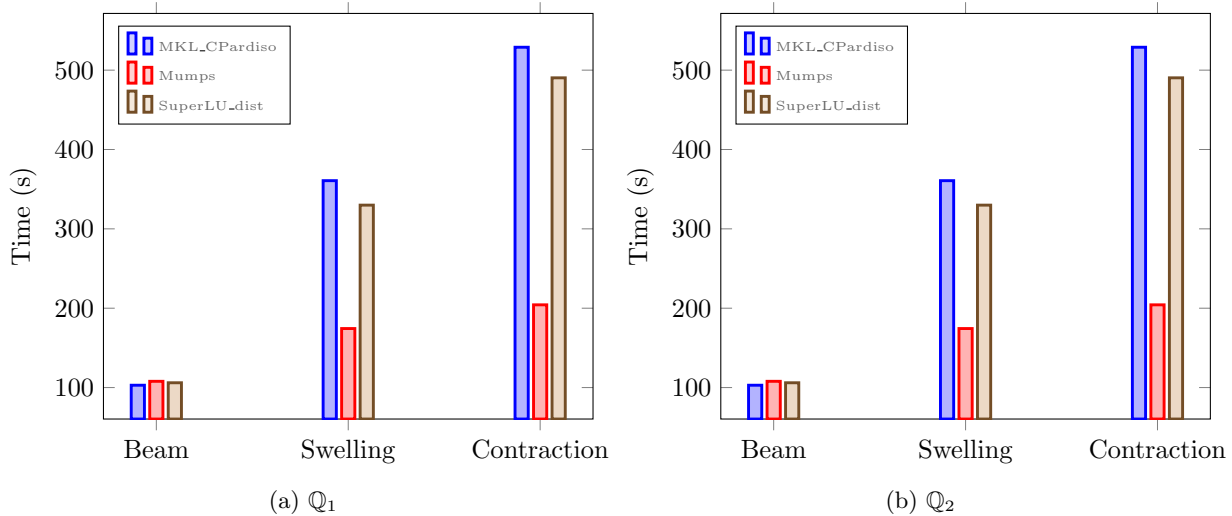
\begin{figure}
    \begin{subfigure}{0.49\textwidth}
        \centering
        \input{tikz/direct/coarse_solver_Q1}
        \caption{\Qone}
    \end{subfigure}
    \begin{subfigure}{0.49\textwidth}
        \centering
        \input{tikz/direct/coarse_solver_Q2}
        \caption{\Qtwo}
    \end{subfigure}
    \caption{CPU time for varying coarse solvers using 32 cores with first and second order finite element spaces.}
    \label{fig:tuning:coarse}
\end{figure}

\paragraph{Coarse problem solver.} We test the performance of the libraries MKL--Pardiso, Mumps, and SuperLU; all of them consider using Mumps for the local problem. We do this in the same setting used in the local solvers comparison, i.e. 15625 degrees of freedom per processor on average. We show the results in Figure \ref{fig:tuning:coarse}, where as before there is little difference between first and second order finite elements. Instead, the biggest difference can be seen between structured ('beam' test) and unstructured ('swelling' and 'contraction' tests). In all cases, the best performance is given by the Mumps library, and in contrast to the serial case, the worst performance is consistently given by MKL--Pardiso.

We conclude the direct solvers study by mentioning that all solvers were tested using an LU decomposition to obtain a fair comparison, but some of them support also a Cholesky factorization. Indeed, only Mumps supports a parallel Cholesky factorization, which further reduces the memory footprint and the CPU time, which makes Mumps the best suited library for using with the BDDC preconditioner in cardiac mechanics in all cases. 

We remark that it could be possible to use inexact or iterative solvers for the local and coarse problems. In fact, in all of our preliminary tests, iterative solvers for the local problems exhibited a significantly worse performance than the direct ones, so we did not consider them for these tests. In the coarse solver, this issue is more delicate. Indeed, as we show in Section \ref{section:scalability}, there is a point in which the communication required for the direct solution of the coarse problem dominates and deteriorates the solution time. This can be alleviated by using multilevel strategies, where the coarse problem is again solved with a BDDC preconditioner computed on a subset of the processors used for the original problem, and allow for extreme scale computing \cite{zampini2016pcbddc}.

\subsection{BDDC primal degrees of freedom}
In this section, we test the performance of different configurations of primal degrees of freedom in the BDDC preconditioner. In particular, we test a configuration where the primal constraints are associated with subdomain vertices (V), subdomain vertices and edges (VE),  subdomain edges and faces (EF), and subdomain vertices, edges and faces (VEF). We note that edges and faces primal constraints are given in terms of averages over the subdomains edges and faces, hence their contribution to the size of the coarse solver is relatively small. This makes it easy to see why the iterations decrease when going from V to VE to VEF primal spaces. Instead, EF shows sometimes improved iteration counts even though there are no vertex constraints, which make up most of the coarse degrees of freedom. All tests are run with 32 cores and roughly one million degrees of freedom. 

\begin{figure}
    \begin{subfigure}{0.49\textwidth}
        \input{tikz/primal/primal_beam_Q1}
        \caption{\Qone}
    \end{subfigure}
    \begin{subfigure}{0.49\textwidth}
        \input{tikz/primal/primal_beam_Q2}
        \caption{\Qtwo}
    \end{subfigure}
    
    \begin{subfigure}{0.49\textwidth}
        \input{tikz/primal/primal_swelling_Q1}
        \caption{\Qone}
    \end{subfigure}
    \begin{subfigure}{0.49\textwidth}
        \input{tikz/primal/primal_swelling_Q2}
        \caption{\Qtwo}
    \end{subfigure}
    
    \begin{subfigure}{0.49\textwidth}
        \input{tikz/primal/primal_realistic_Q1}
        \caption{\Qone}
    \end{subfigure}
    \begin{subfigure}{0.49\textwidth}
        \input{tikz/primal/primal_realistic_Q2}
        \caption{\Qtwo}
    \end{subfigure}
    \caption{Average linear iterations required at each time step in the beam (first row), swelling (second row) and contraction (third row) tests. The legend reports the BDDC primal space: 'V' vertices, 'VE' vertices and edges, 'EF' edges and faces; and 'VEF' vertices, edges and faces. }
    \label{fig:primal}
\end{figure}
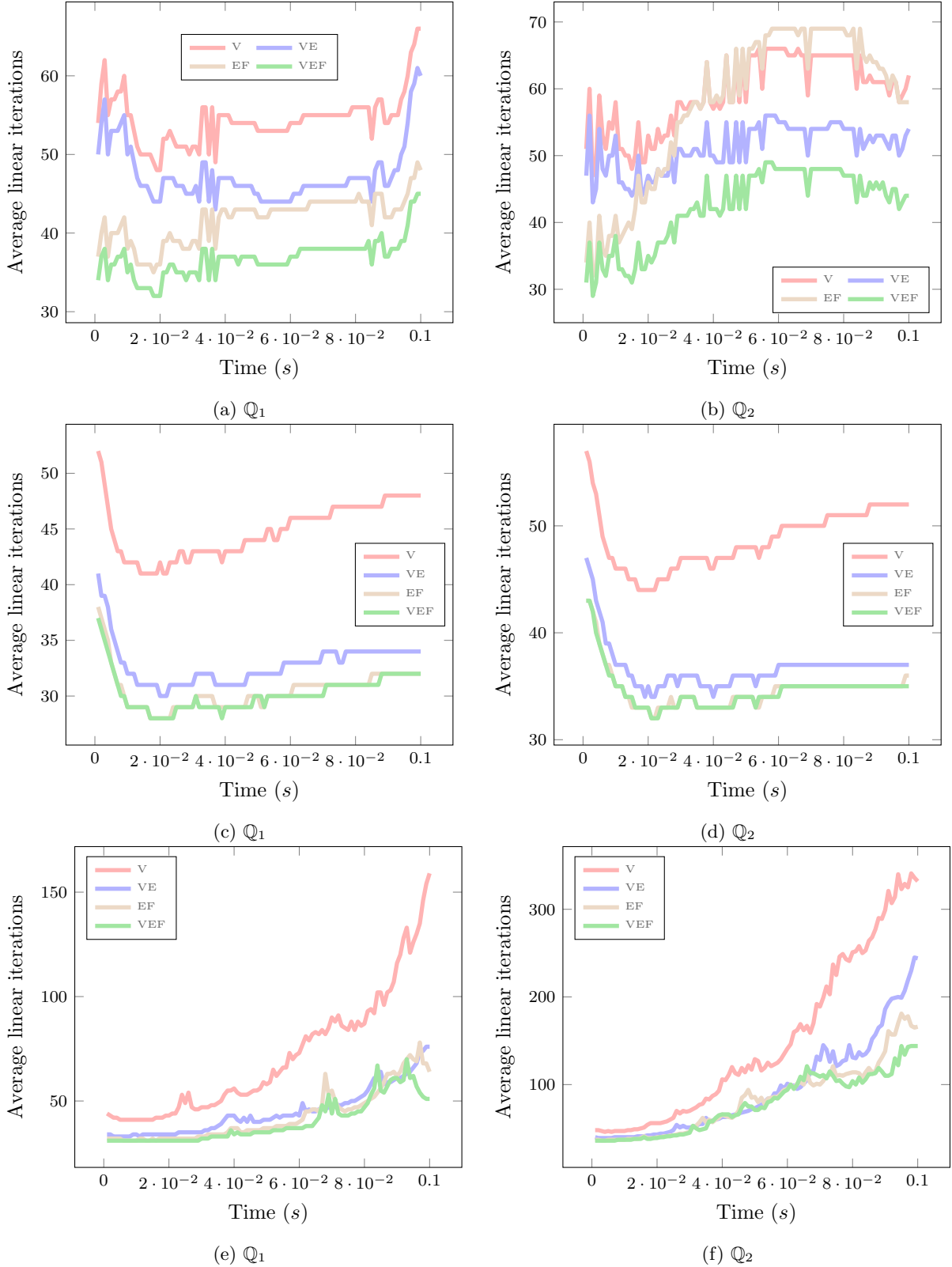

We show the average GMRES iterations throughout the simulation for the beam, swelling and contraction tests in Figure \ref{fig:primal}. We note that, as expected, the richer is the coarse space, the smaller are the iteration counts, with the V, VE and VEF coarse spaces. Instead, interestingly, with the EF space deteriorates over time in the beam test, whereas it outperforms the VE space in the swelling and contraction tests. We also highlight that there is only a small increase in the iterations when going from first to second order finite elements, which conveys the robustness of the BDDC preconditioner when using higher order elements. 
We note that the impact of the size of the coarse space in the overall solution time is not too significant, so that there is no big impact in performance with respect to this choice as shown in \cite{pavarinoSZ2015} (possible reductions in time of up to 30\% in structured meshes), but instead the main aspect to consider here is the robustness and scalability of the preconditioner. Additionally, it is important to weigh the increase in size of the coarse problem with respect to the overall matrix dimension. This impacts the performance when larger problems are considered that require a higher number of cores to be solved. In this case, the coarse problem, which is small, will be distributed among too many processors and eventually communication starts dominating the solution times. The remedy for this are multilevel schemes, more details in Section \ref{section:scalability}.

\section{Strong scalability tests}\label{section:scalability}
In this section we study the strong scalability of our solver in two settings, one with roughly 1 million degrees of freedom (fine mesh) with 2 to 128 CPUs, and one with roughly 8 million degrees of freedom (finer mesh) with 32 to 1024 CPUs. 

In the finer mesh, we have observed a dominant CPU time in the coarse solver of the BDDC when using more than 256 cores (less than 30000 DoFs on average per core), so we also report the CPU time of a multilevel approach (see Section \ref{section:bddc}), where the coarse problem itself is approximated by the application of a BDDC preconditioner on an MPI sub-communicator with the total number of cores divided by 8. For the 1M mesh, we use use edges and faces for the primal space, and instead for the multi-level formulation we use all modes (vertices, edges and faces) to keep the iteration count low.

For all tests, we consider the total CPU time and the total linear iterations during the first 3 time steps.

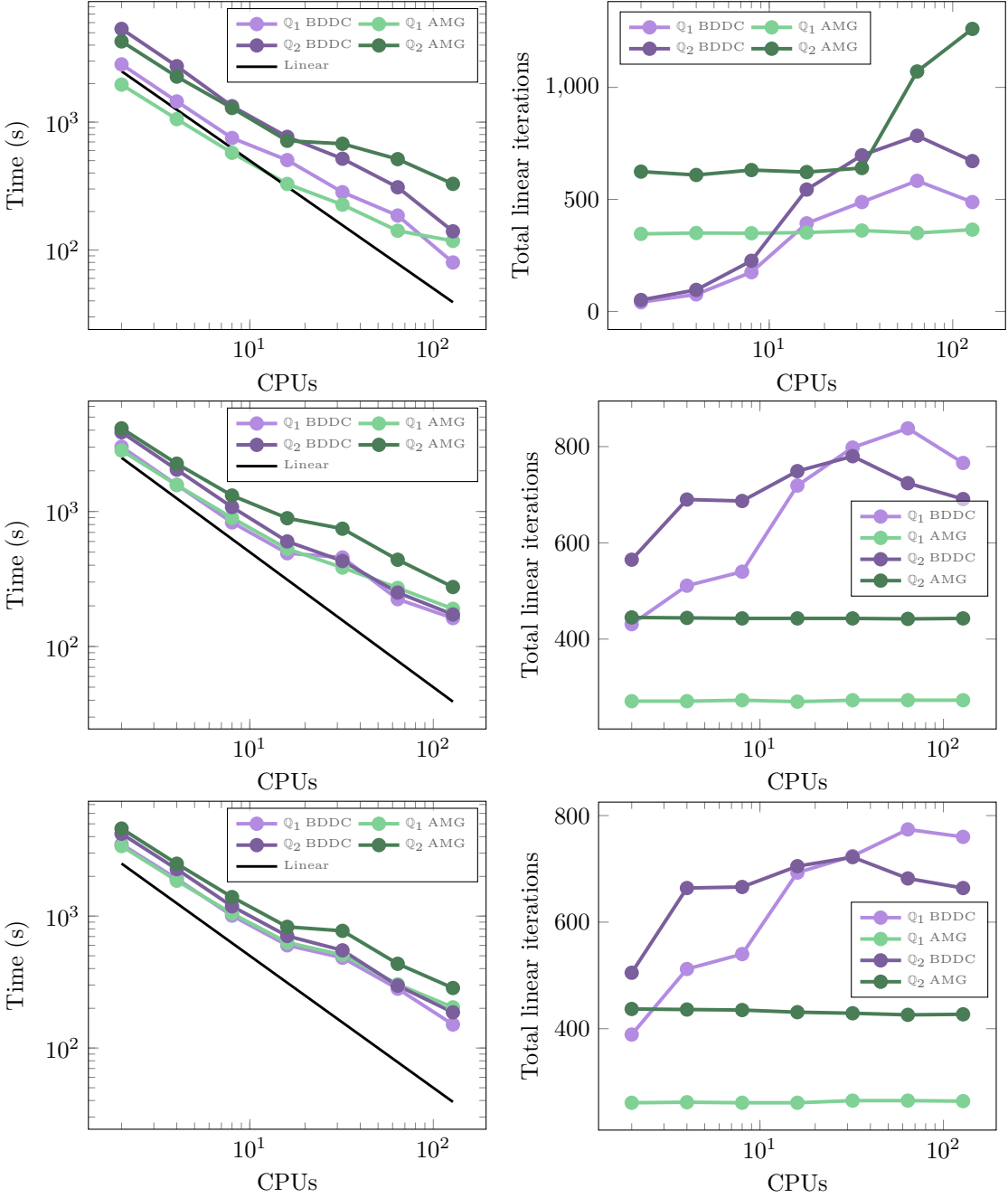
\begin{figure}
    \begin{subfigure}{0.45\textwidth}
    \centering
    \input{tikz/strong/strong_scaling_beam_time}
    \end{subfigure}
    \begin{subfigure}{0.45\textwidth}
    \centering
    \input{tikz/strong/strong_scaling_beam_iterations}
    \end{subfigure}
    
    \begin{subfigure}{0.45\textwidth}
    \centering
    \input{tikz/strong/strong_scaling_swelling_time}
    \end{subfigure}
    \begin{subfigure}{0.45\textwidth}
    \centering
    \input{tikz/strong/strong_scaling_swelling_iterations}
    \end{subfigure}
    
    \begin{subfigure}{0.45\textwidth}
    \centering
    \input{tikz/strong/strong_scaling_contraction_time}
    \end{subfigure}
    \begin{subfigure}{0.45\textwidth}
    \centering
    \input{tikz/strong/strong_scaling_contraction_iterations}
    \end{subfigure}
    \caption{Strong scaling for 1M DoFs, CPU time and total linear iterations in the beam (first row), swelling (second row) and contraction (third row) tests.}
    \label{fig:scaling strong}
\end{figure}

\paragraph{Strong scalability, 1M DoFs.} We show the results for the coarser mesh. The beam test case yields the results shown in Figure \ref{fig:scaling strong} (first row), where it can bee seen that the performance of * degrades with an increasing number of processors, the effect being much more notorious with second order finite elements. Both the BDDC and the AMG (only for first order) are robust with respect to the number of CPUs, but the solution times obtained with the BDDC solver present better scaling, and indeed when using 64 or more processors (32 or more for second order) we see improvements in the solution times when using the BDDC preconditioner. 

The results of the swelling test are shown in Figure \ref{fig:scaling strong} (second row), where in this case all methods yield robust iteration counts with respect to the number of CPUs. The performance of both AMG and BDDC is roughly equal when using first order elements, whereas for second order elements BDDC scales adequately and AMG deteriorates. The results from the contraction test, shown in Figure \ref{fig:scaling strong} (third row) yield exactly the same conclusions, with only the difference that first order is faster with BDDC already with 16 cores, and second order with 8.  

 Note also that in all cases, the number of iterations obtained using the BDDC preconditioner tend to decrease after a certain threshold. This is to be expected, and can be seen from the classic poly-logarithmic bound arising from domain decomposition methods, which depends on $H/h$, i.e. the ration between the largest subdomain size and the discretization size. As the number of cores increases--and so do the subdomains as well-- $h$ remains fixed and $H$ decreases. The bounds are not sharp, so this trend can be seen only after a sufficient number of subdomains is used. We also highlight that the optimized parameters we computed render the preconditioner highly competitive, as can be seen from the ratio between the CPU time of the BDDC and the AMG solvers shown in Table \ref{tab:strong ratios}. Despite a small deterioration in the performance in the unstructured case, there is always  a clear advantage in using the BDDC.

\begin{table}
    \centering
    \begin{tabular}{r | c | c}
        \toprule & \multicolumn{2}{c}{Solution time ratio}\\
        \cmidrule{2-3} & \phantom{asdf}\Qone\phantom{asdf} & \phantom{asdf}\Qtwo\phantom{asdf} \\
        \midrule Beam        & 68\% & 43\% \\
        Swelling    & 85\% & 63\% \\
        Contraction & 74\% & 65\% \\ \bottomrule
    \end{tabular}
    
    \caption{Ratio between the solution times of the BDDC preconditioner and the AMG in the 1M scaling test using 128 cores.}
    \label{tab:strong ratios}
\end{table}

\paragraph{Strong scalability, 8M DoFs.} We show the results for the finer mesh. In this case, we include the results of the multilevel formulation (tagged as BDDC-ML), so we provide the results for Q1 and Q2 separately. The results of the beam test are shown in Figure \ref{fig:scaling stronger beam}. In this case, we note that for Q1 elements, BDDC performs better than AMG only when using from 128 to 512 cores (up to roughly 15625 DoFs per core), then it deteriorates as the solution time of the coarse solver dominates. This can be alleviated using the multilevel method, that presents a satisfactory scaling for all values of CPU cores. This is even clearer when using second order finite elements, where 64 cores is the threshold upon which BDDC is better than AMG. We highlight that the linear iterations of both AMG and BDDC are robust with respect to the CPU cores, but the best scalability is obtained with the multilevel preconditioner, whose iterations deteriorate. In this case, it can be appreciated how AMG is not adequate for second order finite elements. 

\begin{figure}
    \begin{subfigure}{0.45\textwidth}
    \centering
    \input{tikz/stronger/stronger_scaling_swelling_time_Q1}
    \caption{\Qone  CPU time}
    \end{subfigure}
    \begin{subfigure}{0.45\textwidth}
    \centering
    \input{tikz/stronger/stronger_scaling_swelling_iterations_Q1}
    \caption{\Qone  total linear iterations}
    \end{subfigure}
    
    \vspace{0.5cm}
    \begin{subfigure}{0.45\textwidth}
    \centering
    \input{tikz/stronger/stronger_scaling_swelling_time_Q2}
    \caption{\Qtwo  CPU time}
    \end{subfigure}
    \begin{subfigure}{0.45\textwidth}
    \centering
    \input{tikz/stronger/stronger_scaling_swelling_iterations_Q2}
    \caption{\Qtwo  total linear iterations}
    \end{subfigure}
    \caption{Strong scaling for 8M DoFs, swelling test.}
    \label{fig:scaling stronger swelling}
\end{figure}
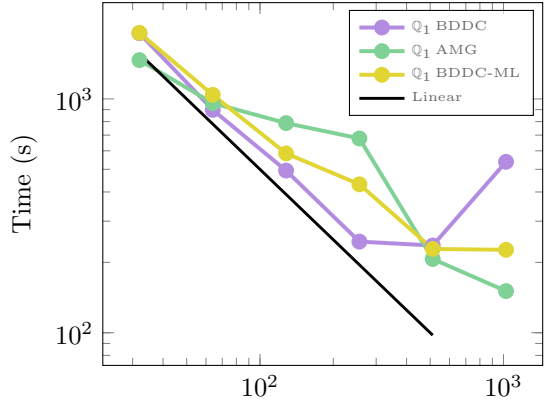
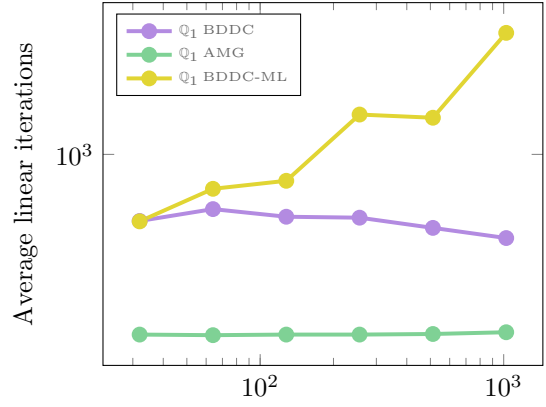
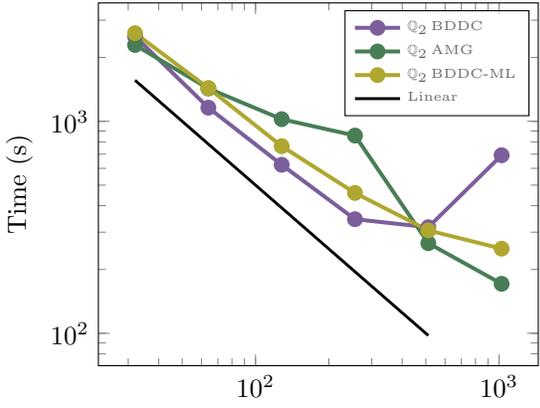
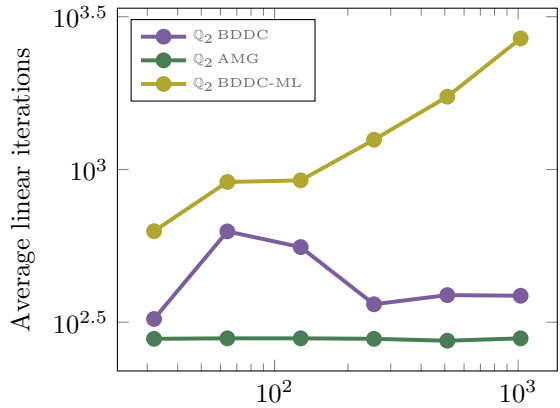

The results for the swelling test are shown in Figure \ref{fig:scaling stronger swelling}. In this test we obtained that again both BDDC and AMG present a robust iteration count with respect to the CPUs and the multilevel instead deteriorates as the CPUs increase. In this test, the BDDC preconditioner scales adequately only up to 256 cores in both Q1 and Q2 formulations (31250 DoFs per core on average). After that, AMG performs better. Note that this difficulty could not be circumvented with the multilevel preconditioner. 

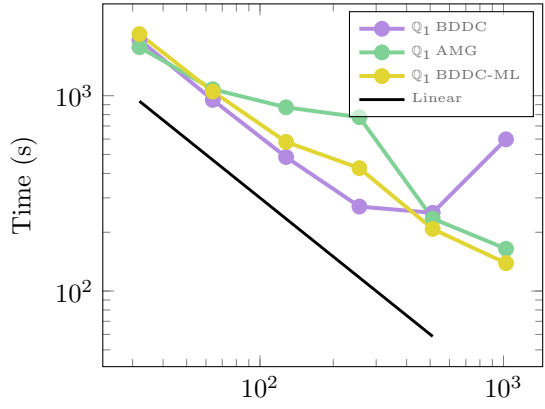
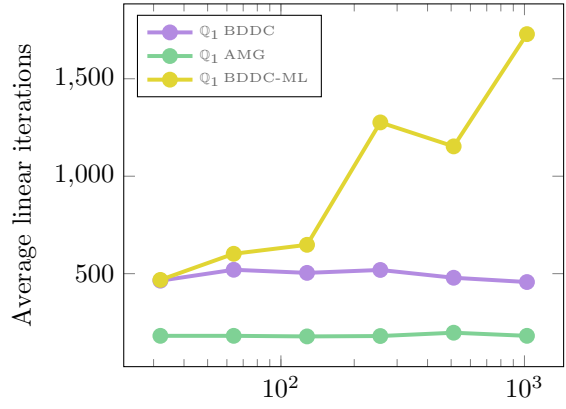
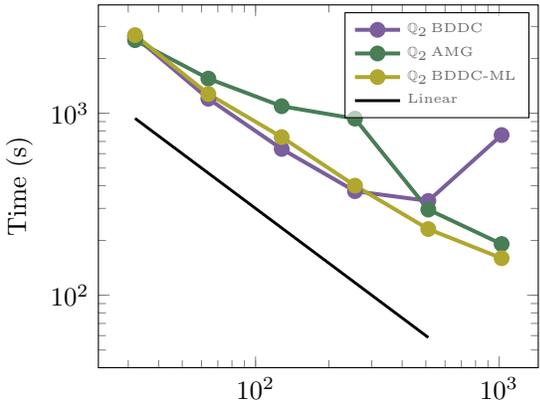
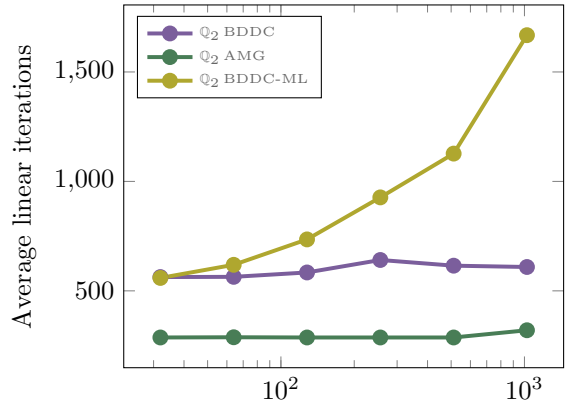
\begin{figure}
    \begin{subfigure}{0.45\textwidth}
    \centering
    \input{tikz/stronger/stronger_scaling_contraction_time_Q1}
    \caption{\Qone  CPU time}
    \end{subfigure}
    \begin{subfigure}{0.45\textwidth}
    \centering
    \input{tikz/stronger/stronger_scaling_contraction_iterations_Q1}
    \caption{\Qone  total linear iterations}
    \end{subfigure}
    
    \vspace{0.5cm}
    \begin{subfigure}{0.45\textwidth}
    \centering
    \input{tikz/stronger/stronger_scaling_contraction_time_Q2}
    \caption{\Qtwo  CPU time}
    \end{subfigure}
    \begin{subfigure}{0.45\textwidth}
    \centering
    \input{tikz/stronger/stronger_scaling_contraction_iterations_Q2}
    \caption{\Qtwo  total linear iterations}
    \end{subfigure}
    \caption{Strong scaling for 8M DoFs, contraction test.}
    \label{fig:scaling stronger contraction}
\end{figure}

Finally, the results for the contraction test are shown in Figure \ref{fig:scaling stronger contraction}. In this case, the qualitative results obtained in the swelling case regarding the iteration counts are identical, but instead solution times are different. Indeed, it is still true that pure BDDC presents the best scalability with up to 256 cores, but after that the multilevel preconditioner presents the best performance. 

In this case, the percentages of CPU time of the BDDC with respect to the AMG preconditioner must be considered in two separate cases, as plain BDDC only scales up to 256 cores in the unstructured case. So, we compare the performance of the BDDC preconditioner in the 256 cores scenario, and the multilevel in the 1024 cores scenario as shown in Table \ref{tab:stronger ratios}. We highlight that in roughly all cases, the BDDC preconditioner performs better than AMG. Still, for the 1024 cores case, AMG performs better than BDDC-ML only in the swelling case. The advantage of using a multilevel BDDC method is less striking than the one obtained using plain BDDC, at least in the ranges where it performs better. This can be already explained by the deterioration in the iteration count yielded by the BDDC-ML, and indeed improving the coarse solver is one of our main future objectives.

\begin{table}
    \centering
    \begin{subtable}{0.49\textwidth}
        \begin{tabular}{r | c | c}
        \toprule & BDDC & BDDC-ML  \\ \midrule
        Beam        & 44\%  & 72\% \\
        Swelling    & 36\%  & 150\% \\
        Contraction & 40\%  & 84\%  \\ \bottomrule
        \end{tabular}
        \caption{$\mathbb Q_1$}
    \end{subtable}
    \begin{subtable}{0.49\textwidth}
        \begin{tabular}{r | c | c}
        \toprule & BDDC & BDDC-ML \\ \midrule
        Beam        & 26\%  & 20\% \\
        Swelling    & 40\%  & 147\% \\
        Contraction & 35\%  & 84\% \\ \bottomrule
        \end{tabular}
        \caption{$\mathbb Q_2$}
    \end{subtable}
    \caption{Ratio of the solution time with the BDDC preconditioner with respect to the AMG in the 8M scaling test. In plain BDDC, the coarse problem dominates the solution time when using more than 256 cores, so we display the results for that scenario, and instead compare the results using the multilevel method (BDDC-ML) in the 1024 cores scenario.}
    \label{tab:stronger ratios}
\end{table}

\begin{figure}
    \begin{subfigure}{0.45\textwidth}
    \centering
    \input{tikz/stronger/stronger_scaling_beam_time_Q1}
    \caption{\Qone  CPU time}
    \end{subfigure}
    \begin{subfigure}{0.45\textwidth}
    \centering
    \input{tikz/stronger/stronger_scaling_beam_iterations_Q1}
    \caption{\Qone  total linear iterations}
    \end{subfigure}
    
    \vspace{0.5cm} 
    \begin{subfigure}{0.45\textwidth}
    \centering
    \input{tikz/stronger/stronger_scaling_beam_time_Q2}
    \caption{\Qtwo  CPU time}
    \end{subfigure}
    \begin{subfigure}{0.45\textwidth}
    \centering
    \input{tikz/stronger/stronger_scaling_beam_iterations_Q2}
    \caption{\Qtwo  total linear iterations}
    \end{subfigure}
    \caption{Strong scaling for 8M DoFs, beam test.}
    \label{fig:scaling stronger beam}
\end{figure}
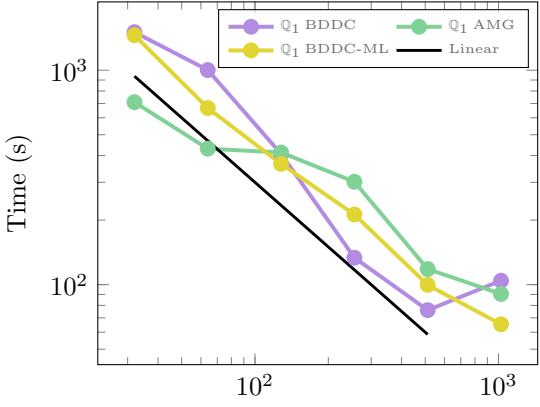
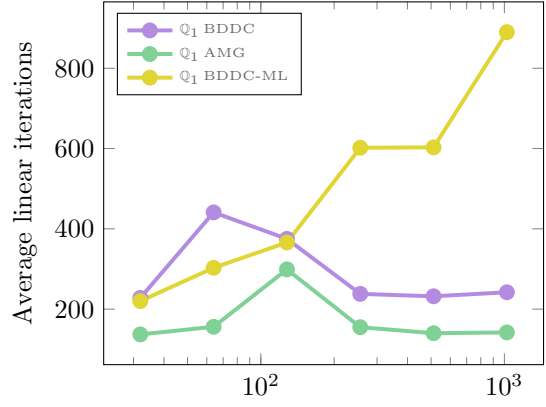
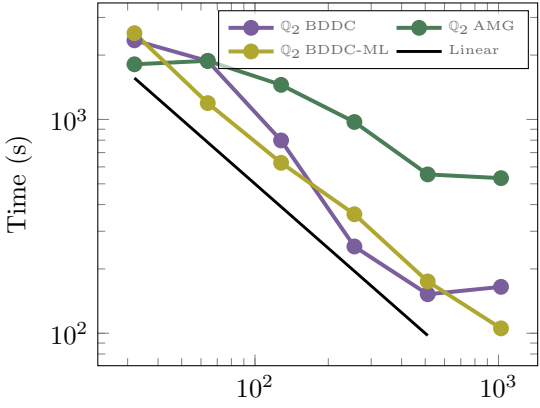
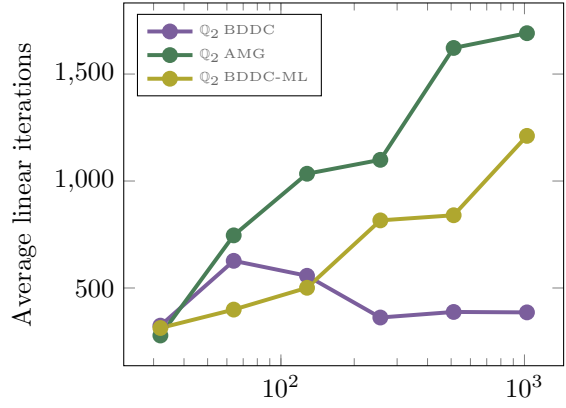

\section{A realistic simulation}\label{section:realistic}
In this test we present a realistic simulation given by the model detailed in \cite{regazzoni2020cardiac}. In addition to the details shown in Section \ref{section:mechanics}, this model considers:
\begin{itemize}
    \item A reduced, closed loop circulation model that describes blood flow through the body and the valve dynamics, which yields a realistic PV loop.
    \item A Holzapfel-Ogden potential instead of the Guccione one.
    \item An activation function given by the evolution of an electrophysiology model. 
\end{itemize}
We refer the interested reader to the original reference \cite{regazzoni2020cardiac} for the details of the different components of the model and the description of numerical methods that are well suited for their solution. In this section we consider only the evolution of one heartbeat, and assess the performance of the BDDC and AMG preconditioners for both first and second order finite elements, to see the difference in the performance. Due to limitations in our server infrastructure, we were not able to perform an entire electromechanical simulation with a BDDC preconditioner with a million DoFs. For this reason, we have devised two scenarios: one small problem ($10^5$ DoFs) and one large ($8\cdot 10^5$ DoFs), both with first order finite elements. In all cases, we consider the BDDC preconditioner with all primal modes (vertices, edges and faces) with an LU decomposition as direct solver for both local and coarse problems.

We remark that this test is biased towards the AMG preconditioner. This happens because the model under consideration uses an average pressure at the base $\partial\Omega_\text{base}$ in order to balance the pressure acting on the endocardium $\partial \Omega_\text{base}$, which ultimately breaks the symmetry of the problem (see \cite{regazzoni2020cardiac} for details on this boundary condition). This does not allow us to use a Cholesky factorization in the local and coarse problems as done in the scalability tests, as well as requiring to inform the BDDC preconditioner to use the non-symmetric formulation, which is more expensive to setup and to apply. In the future, we are interested in performing similar studies on a four-chamber cardiac mechanics setting, in which we can again exploit the symmetry of the mechanics problem.

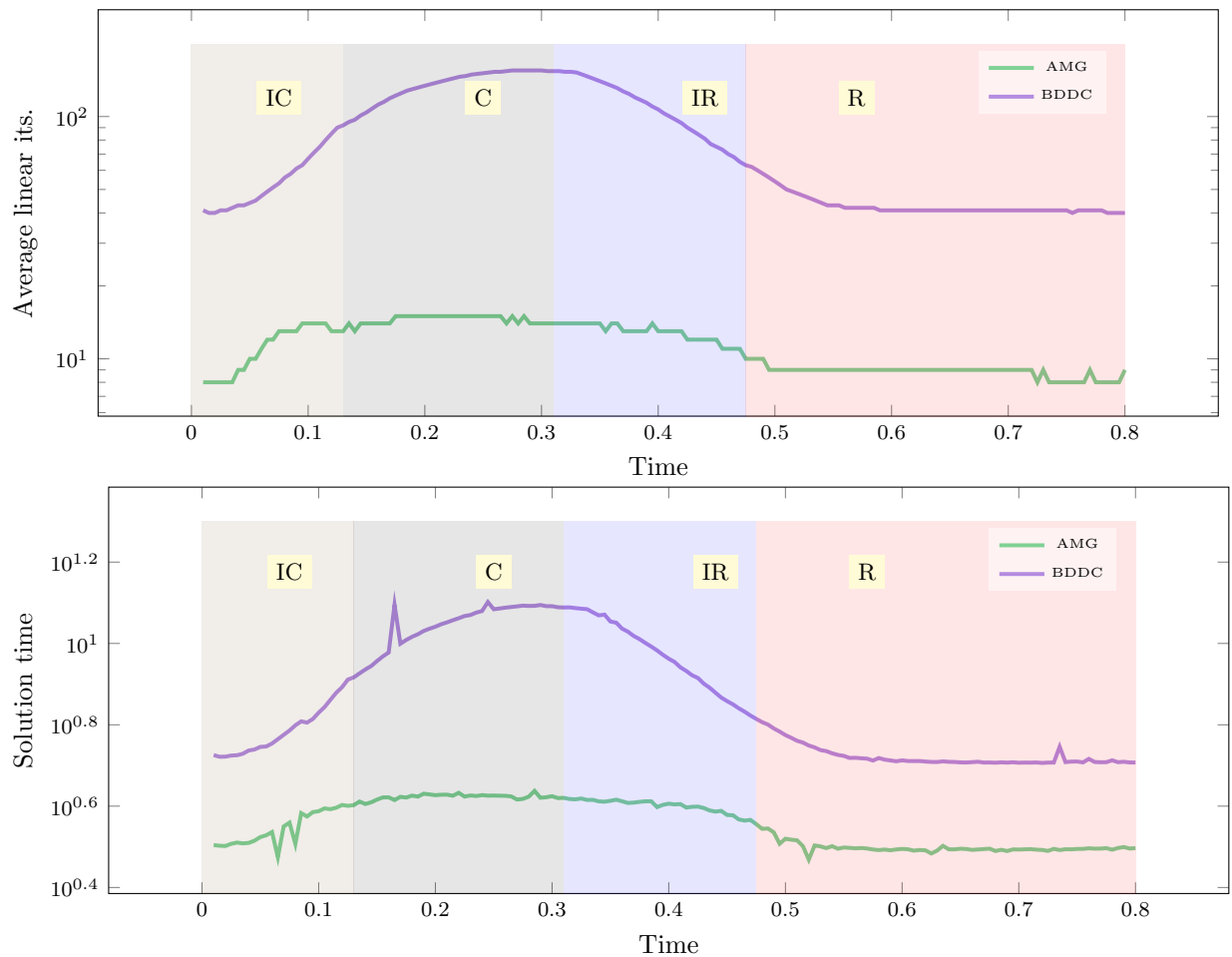
\begin{figure}
    \centering
    \input{tikz/EM/small-p1.tex}
    \caption{Simulation of an entire electromechanical simulation using both an AMG and a BDDC preconditioner for the mechanics. The top image shows the average linear iterations incurred during the Newton steps, the bottom image shows the total time it took to solve the given timestep.}
    \label{fig:heartbeat-small-p1}
\end{figure}

\begin{figure}
    \centering
    \input{tikz/EM/large-p1.tex}
    \caption{Simulation of an entire electromechanical simulation using both an AMG and a BDDC preconditioner for the mechanics with second order elements. The left image shows the average linear iterations incurred during the Newton steps, the right image shows the total time it took to solve the given timestep.}
    \label{fig:heartbeat-large-p1}
\end{figure}
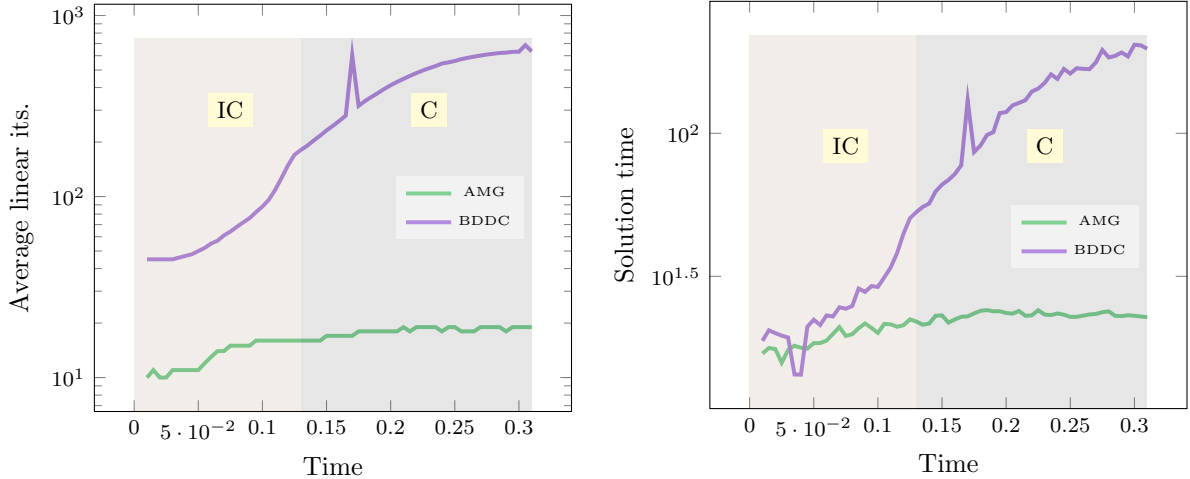

\paragraph{Small simulation.} In this scenario we use 48 CPU cores, with roughly $10^5$ degrees of freedom, the results are reported in Figure \ref{fig:heartbeat-small-p1}, where we have separated the results into the four phases of the heartbeat, i.e. isovolumic contraction (IC), contraction (C), isovolumic relaxation (IR) and final relaxation (R). We note that both methods follow the same trend, meaning that they perform better during the IC and R phases, and instead deteriorate their performance during the C and IR phases. AMG presents an increase in the iteration count of roughly double, going from 11 iterations per time step in average up to 20 in average. The deterioration of BDDC is worse, going from 44 to over 400 iterations in average per time step. This effect greatly impacts the solution time, as can bee seen from the bottom figure in Figure \ref{fig:heartbeat-small-p1}.
\paragraph{Large simulation.} In this scenario we use 120 CPU cores, with roughly $8\cdot 10^5$ degrees of freedom. This test present a similar trend to the one shown in the small simulation, as shown in Figure \ref{fig:heartbeat-large-p1}. In this case, the difference in CPU time between these two methods is smaller during the IC phase, and instead towards the C phase the deterioration of the BDDC preconditioner greatly impacts the solution time, making it less competitive. The improvement of the BDDC preconditioner with respect to the AMG, particularly during the IC phase, can be better understood by its superior scalability, as shown in the benchmark tests in Section \ref{section:scalability}. Still, AMG vastly outperforms BDDC throughout the entire performance, with only a small exception during the IC phase.


\section{Conclusions}

In this work we have performed a detailed comparison study between BDDC and AMG preconditioners in the context of cardiac mechanics. We have observed a great impact in the overall solution time using the BDDC preconditioner due to the choice of the direct solver used in both local and coarse problems, with MUMPS being the best library to use in this context. We have also observed that in some cases it can show little to no difference to use the vertices in the primal space for BDDC, as observed in the benchmark tests. 

The performance of the BDDC preconditioner in the benchmark tests is more than adequate, and in general it outperforms the AMG. Still, this is not the case for the study of a complete heartbeat, where AMG proves to be much more robust. Both preconditioners perform worse during the contraction (C) and isovolumic relaxation (IR) phases, but the deterioration in performance of BDDC is much worse. Still, it needs to be further studied how symmetry impacts the performace of the solver, as scalability was performed in symmetric problems that allowed for the use of a Cholesky factorization in BDDC, whereas the heartbeat simulation allows only for a generic LU factorization. 

\section*{Acknowledgments}
We would like to thank Stefano Zampini for his help in making an efficient use of the BDDC preconditioner in PETSc. N. Barnafi and L. F. Pavarino have been supported by
grants of MIUR (PRIN 2017AXL54F$\_$002) and INdAM--GNCS. 
N. Barnafi and S. Scacchi have been supported by grants of MIUR (PRIN
2017AXL54F$\_$003) and INdAM-GNCS. 
The Authors are also grateful to the University of Pavia, the University of Milan, and the CINECA laboratory for the usage of the EOS, INDACO and Galileo100 clusters, respectively
\bibliography{main}
\bibliographystyle{alpha}
\end{document}

%% file: preamble.tex
\usepackage[utf8]{inputenc}
\usepackage{amsmath, amssymb}
\usepackage{mathtools} 
\usepackage{fullpage}
\usepackage{todonotes}
\usepackage{subcaption} 
\usepackage{booktabs}
\usepackage{isomath}
\usepackage{pgfplots}  
\usepackage{dsfont} 
\pgfplotsset{compat=1.15}
\pgfplotsset{subindexes/.style 2 args={
    x filter/.code={
        \ifnum\coordindex<#1\fi
        \ifnum\coordindex>#2\fi
    }
}}
\pgfplotsset{my line/.style={%
        line width=2pt,line join=round,color=#1, subindexes={1}{100}
    }
}

\renewcommand*{\vec}{\vectorsym}
\newcommand{\ten}{\tensorsym}
\newcommand{\mat}{\matrixsym}

\newcommand{\R}{\mathbb{R}}
\newcommand{\Q}{\mathbb{Q}}
\newcommand{\Qone}{\ensuremath{\Q_{1}\,}}
\newcommand{\Qtwo}{\ensuremath{\Q_{2}\,}}

\newcommand{\parder}[2]{\frac{\partial\,#1}{\partial\,#2}}
\DeclareMathOperator{\dive}{\text{div}}

\newcommand{\bbone}{\mathds{1}}

\newcommand{\displacement}{\vec{d}}
\newcommand{\tenF}{\ten{F}}

\pgfplotsset{every axis plot/.append style={line width=1.1pt}}

\newenvironment{semilogyplot}[3][]
      {\begin{tikzpicture}
      \begin{semilogyaxis}[width=\textwidth, xlabel=#2, ylabel=#3, tick label style={font=\footnotesize}, legend style={font=\tiny, draw=none, fill opacity=0.5, text opacity = 1,row sep=2pt},legend pos=north west, legend columns=1, #1]
      }
      {%
      \end{semilogyaxis}
      \end{tikzpicture}
      }

\newcommand{\codecolor}{\color{purple!80!white}}
\newcommand{\icode}[1]{{\codecolor \texttt{#1}}} 

\usepackage{listings}
\usepackage{xcolor}
\definecolor{codegreen}{rgb}{0,0.6,0}
\definecolor{codegray}{rgb}{0.5,0.5,0.5}
\definecolor{codepurple}{rgb}{0.58,0,0.82}
\definecolor{backcolour}{rgb}{0.95,0.95,0.92}
\lstdefinestyle{mystyle}{
  backgroundcolor=\color{backcolour}, commentstyle=\color{codegreen},
  keywordstyle=\color{magenta},
  numberstyle=\tiny\color{codegray},
  stringstyle=\color{codepurple},
  basicstyle=\ttfamily\small,
  breakatwhitespace=false,         
  captionpos=b,                    
  keepspaces=true,                 
  numbers=none,                    
  numbersep=5pt,                  
  showspaces=false,                
  showstringspaces=false,
  showtabs=false,                  
  tabsize=2,
  framerule=1.5pt,
  rulecolor=\color{red!60!black}
}

\pgfplotsset{select coords between index/.style 2 args={
    x filter/.code={
        \ifnum\coordindex<#1\fi
        \ifnum\coordindex>#2\fi
    }
}}

\definecolor{bddcp1}{RGB}{179, 139, 225}
\definecolor{bddcp2}{RGB}{123, 94, 156}
\definecolor{bddcmlp1}{RGB}{225, 213, 51}
\definecolor{bddcmlp2}{RGB}{175, 167, 47}
\definecolor{amgp1}{RGB}{127, 209, 149}
\definecolor{amgp2}{RGB}{72, 125, 86}
\pgfplotsset{BDDCP1style/.style={
        mark options={draw opacity=0.0, fill opacity=1.0, mark size=3pt}, 
        mark=*,
        color=bddcp1, 
        line width=1.5pt}
        }
\pgfplotsset{BDDCP2style/.style={
        mark options={draw opacity=0.0, fill opacity=1.0, mark size=3pt}, 
        mark=*,
        color=bddcp2, 
        line width=1.5pt}
        }
\pgfplotsset{AMGP1style/.style={
        mark options={draw opacity=0.0, fill opacity=1.0, mark size=3pt}, 
        mark=*,
        color=amgp1, 
        line width=1.5pt}
        }
\pgfplotsset{AMGP2style/.style={
        mark options={draw opacity=0.0, fill opacity=1.0, mark size=3pt}, 
        mark=*,
        color=amgp2, 
        line width=1.5pt}
        }
\pgfplotsset{BDDCMLP1style/.style={
        mark options={draw opacity=0.0, fill opacity=1.0, mark size=3pt}, 
        mark=*,
        color=bddcmlp1, 
        line width=1.5pt}
        }
\pgfplotsset{BDDCMLP2style/.style={
        mark options={draw opacity=0.0, fill opacity=1.0, mark size=3pt}, 
        mark=*,
        color=bddcmlp2, 
        line width=1.5pt}
        }

%% file: arxiv-header.tex
\title{A comparative study of scalable multilevel preconditioners for cardiac mechanics}
\author{Nicolás A. Barnafi\thanks{Department of Mathematics, Universit\`a di Pavia, Via Ferrata 1, 27100 Pavia, Italy. \texttt{\{nicolas.barnafi, luca.pavarino\}@unipv.it}} \and Luca F. Pavarino\footnotemark[1] \and Simone Scacchi\thanks{Department of Mathematics, Universit\`a di Milano,
Via Saldini 50, 20133 Milano, Italy. \texttt{simone.scacchi@unimi.it}}}
\date{}
\maketitle
\begin{abstract}
    In this work, we provide a performance comparison between the Balancing Domain Decomposition by Constraints (BDDC) and the Algebraic Multigrid (AMG) preconditioners for cardiac mechanics on both structured and unstructured finite element meshes. The mechanical behavior of myocardium can be described by the equations of three-dimensional finite elasticity, which are discretized by finite elements in space and yield the solution of a large scale nonlinear algebraic system. This problem is solved by a Newton-Krylov method, where the solution of the Jacobian linear system is accelerated by BDDC/AMG preconditioners. We thoroughly explore the main parameters of the BDDC preconditioner in order to make the comparison fair. We focus on: the performance of different direct solvers for the local and coarse problems of the BDDC algorithm; the impact of the different choices of BDDC primal degrees of freedom; and the influence of the finite element degree. Scalability tests are performed on Linux clusters up to 1024 processors, and we conclude with a performance study on a realistic electromechanical simulation.
\end{abstract}

%% file: els-header.tex
\begin{frontmatter}
\title{A comparative study of scalable multilevel preconditioners for cardiac mechanics}

\author[1]{Nicolás A. Barnafi}
\author[1]{Luca F. Pavarino}
\author[2]{Simone Scacchi}

\address[1]{Department of Mathematics, Universit\`a degli Studi di Pavia, Italy. \texttt{\{nicolas.barnafi,luca.pavarino\}@unipv.it}}
\address[2]{Department of Mathematics, Universit\`a degli Studi di Milano, Italy. \texttt{simone.scacchi@unimi.it}}

\begin{abstract}
In this work, we provide a performance comparison between the Balancing Domain Decomposition by Constraints (BDDC) and the Algebraic Multigrid (AMG) preconditioners for cardiac mechanics on both structured and unstructured finite element meshes. The mechanical behavior of myocardium can be described by the equations of three-dimensional finite elasticity, which are discretized by finite elements in space and yield the solution of a large scale nonlinear algebraic system. This problem is solved by a Newton-Krylov method, where the solution of the Jacobian linear system is accelerated by BDDC/AMG preconditioners. We thoroughly explore the main parameters of the BDDC preconditioner in order to make the comparison fair. We focus on: the performance of different direct solvers for the local and coarse problems of the BDDC algorithm; the impact of the different choices of BDDC primal degrees of freedom; and the influence of the finite element degree. Scalability tests are performed on Linux clusters up to 1024 processors, and we conclude with a performance study on a realistic electromechanical simulation.
\end{abstract}

\begin{keyword}
\MSC 65F08 \sep 65M22 \sep 65N20 
\end{keyword}



\end{frontmatter}

%% file: tikz/direct/local_solver_parallel_Q1.tex
    \begin{tikzpicture}
    \begin{axis}
        [ybar, xtick=data, width=\textwidth, height=6cm, ymin=0,
        ylabel=Time (s), bar width=5pt, enlarge x limits={abs=20pt},  
        symbolic x coords={Beam, Swelling, Contraction},
         legend style={at={(0.03,0.97)},anchor=north west, fill opacity=0.6, font=\tiny,legend cell align=left}
        ]
        \addplot coordinates {(Beam,699.372) 
                              (Swelling, 800.129)
                              (Contraction, 912.039)};
        \addplot coordinates {(Beam, 715.614) 
                              (Swelling,  729.923)
                              (Contraction, 817.680)};
        \addplot coordinates {(Beam, 1298.235)
                              (Swelling, 1327.576)
                              (Contraction, 1419.990)};
        \addplot coordinates {(Beam, 2610.109)
                              (Swelling, 4955.759)
                              (Contraction, 5103.388)};
                              
        \addplot[fill=yellow!40!white, draw=yellow!50!black] coordinates {(Beam, 6332.620)
                              (Swelling, 3394.630)
                              (Contraction, 3576.072)};
        \addplot[fill=green!40!white, draw=green!50!black] coordinates {(Beam,3886.112)
                              (Swelling, 21570.941)
                              (Contraction, 22367.597)};
    \legend{MKL\_Pardiso, Mumps, UMFPACK, SuperLU, KLU, PETSc LU}
    \end{axis}
    \end{tikzpicture}
    

%% file: tikz/direct/local_solver_parallel_Q2.tex
    \begin{tikzpicture}
    \begin{axis}
        [ybar,
        xtick=data,
        width=\textwidth,
        height=6cm,
        ymin=0,
        ylabel=Time (s),
        bar width=5pt,
        enlarge x limits={abs=20pt},  
        symbolic x coords={Beam, Swelling, Contraction},
         legend style={at={(0.03,0.97)},anchor=north west, fill opacity=0.6, font=\tiny,legend cell align=left}
        ]
        \addplot coordinates {(Beam,1238.733) 
                              (Swelling, 952.455)
                              (Contraction, 1023.837)};
        \addplot coordinates {(Beam, 1240.032) 
                              (Swelling,  899.006)
                              (Contraction, 953.329)};
        \addplot coordinates {(Beam, 1926.608)
                              (Swelling, 1403.091)
                              (Contraction, 1520.020)};
        \addplot coordinates {(Beam,7293.498)
                              (Swelling, 15472.934)
                              (Contraction, 15896.541)};
        \addplot[fill=yellow!40!white, draw=yellow!50!black] coordinates {(Beam, 7677.900)
                              (Swelling, 0.0)
                              (Contraction, 1215.944)};
        \addplot[fill=green!40!white, draw=green!50!black] coordinates {(Beam,6217.853)
                              (Swelling, 18128.609)
                              (Contraction, 19085.064)};
    \legend{MKL\_Pardiso, Mumps, SuperLU, UMFPACK, KLU, PETSc LU}
    \end{axis}
    \end{tikzpicture}

%% file: tikz/direct/local_solver_parallel_Q1_opt.tex
    \begin{tikzpicture}
    \begin{axis}
        [ybar, xtick=data, width=\textwidth, height=6cm, ymin=0,
        ymax=1600,
        ylabel=Time (s), bar width=5pt, enlarge x limits={abs=20pt},  
        symbolic x coords={Beam, Swelling, Contraction},
         legend style={at={(0.03,0.97)},anchor=north west, fill opacity=0.6, font=\tiny,legend cell align=left},
         legend columns=3
        ]
        \addplot coordinates {(Beam,699.372) 
                              (Swelling, 800.129)
                              (Contraction, 912.039)};
        \addplot coordinates {(Beam, 715.614) 
                              (Swelling,  729.923)
                              (Contraction, 817.680)};

        \addplot coordinates {(Beam, 1298.235)
                              (Swelling, 1327.576)
                              (Contraction, 1419.990)};

    \legend{MKL\_Pardiso, Mumps, UMFPACK}
    \end{axis}
    \end{tikzpicture}
    

%% file: tikz/direct/local_solver_parallel_Q2_opt.tex
    \begin{tikzpicture}
    \begin{axis}
        [ybar,
        xtick=data,
        width=\textwidth,
        height=6cm,
        ymin=0,
        ylabel=Time (s),
        bar width=5pt,
        enlarge x limits={abs=20pt},  
        symbolic x coords={Beam, Swelling, Contraction},
         legend style={at={(0.03,0.97)},anchor=north west, fill opacity=0.6, font=\tiny,legend cell align=left}, 
        legend pos=north east, 
        legend columns=2
        ]
        \addplot coordinates {(Beam,1238.733) 
                              (Swelling, 952.455)
                              (Contraction, 1023.837)};
        \addplot coordinates {(Beam, 1240.032) 
                              (Swelling,  899.006)
                              (Contraction, 953.329)};

        \addplot coordinates {(Beam, 1926.608)
                              (Swelling, 1403.091)
                              (Contraction, 1520.020)};

    \legend{MKL\_Pardiso, Mumps, UMFPACK}
    \end{axis}
    \end{tikzpicture}

%% file: tikz/direct/coarse_solver_Q1.tex
    \begin{tikzpicture}
    \begin{axis}
        [ybar,
        xtick=data,
        width=\textwidth,
        ylabel=Time (s),
        bar width=5pt,
        enlarge x limits={abs=20pt},  
        symbolic x coords={Beam, Swelling, Contraction},
         legend style={at={(0.03,0.97)},anchor=north west, fill opacity=0.6, font=\tiny,legend cell align=left}
        ]
        \addplot coordinates {(Beam,102.945) 
                              (Swelling, 360.727)
                              (Contraction, 528.861)};
        \addplot coordinates {(Beam, 107.875) 
                              (Swelling, 174.372)
                              (Contraction, 204.254)};
        SuperLU
        \addplot coordinates {(Beam,106.051)
                              (Swelling, 329.892)
                              (Contraction,490.316)};
    \legend{MKL\_CPardiso, Mumps, SuperLU\_dist}
    \end{axis}
    \end{tikzpicture}




%% file: tikz/direct/coarse_solver_Q2.tex
    \begin{tikzpicture}
    \begin{axis}
        [ybar,
        xtick=data,
        width=\textwidth,
        ylabel=Time (s),
        bar width=5pt,
        enlarge x limits={abs=20pt},  
        symbolic x coords={Beam, Swelling, Contraction},
         legend style={at={(0.03,0.97)},anchor=north west, fill opacity=0.6, font=\tiny,legend cell align=left}
        ]
        \addplot coordinates {(Beam,102.945) 
                              (Swelling, 360.727)
                              (Contraction, 528.861)};
        \addplot coordinates {(Beam, 107.875) 
                              (Swelling, 174.372)
                              (Contraction, 204.254)};
        SuperLU
        \addplot coordinates {(Beam,106.051)
                              (Swelling, 329.892)
                              (Contraction,490.316)};
    \legend{MKL\_CPardiso, Mumps, SuperLU\_dist}
    \end{axis}
    \end{tikzpicture}




%% file: tikz/primal/primal_beam_Q1.tex
\newcommand{\indMinBeam}{2}
\newcommand{\indMaxBeam}{100}

\begin{tikzpicture}
\begin{axis}
[width=\textwidth, no markers, xlabel=Time $(s)$, ylabel=Average linear iterations, tick label style={font=\footnotesize}, legend style={fill opacity=0.6, font=\tiny,legend cell align=left, at={(0.5,0.9)}, anchor=north}, legend columns=2]
\addplot [my line=red!30!white] table [x=time, y=mechanics_iterations_linear_avg, col sep=comma] {data/beam/primal-100-Q1.csv};
\addplot [my line=blue!30!white] table [x=time, y=mechanics_iterations_linear_avg, col sep=comma] {data/beam/primal-110-Q1.csv};
\addplot [my line=brown!30!white] table [x=time, y=mechanics_iterations_linear_avg, col sep=comma] {data/beam/primal-011-Q1.csv};
\addplot [my line=green!30!white!90!black] table [x=time, y=mechanics_iterations_linear_avg, col sep=comma] {data/beam/primal-111-Q1.csv};
\legend {V, VE, EF, VEF}
\end{axis}
\end{tikzpicture}

%% file: tikz/primal/primal_beam_Q2.tex
\newcommand{\indMinBeam}{2}
\newcommand{\indMaxBeam}{100}
\pgfplotsset{my line/.style={%
        line width=2pt,line join=round,color=#1, subindexes={1}{100}
    }
}
\begin{tikzpicture}
\begin{axis}
[width=\textwidth, no markers, xlabel=Time $(s)$, ylabel=Average linear iterations, tick label style={font=\footnotesize}, legend style={fill opacity=0.6, font=\tiny,legend cell align=left}, legend pos=south east, legend columns=2]
\addplot [my line=red!30!white] table [x=time, y=mechanics_iterations_linear_avg, col sep=comma] {data/beam/primal-100-Q2.csv};
\addplot [my line=blue!30!white] table [x=time, y=mechanics_iterations_linear_avg, col sep=comma] {data/beam/primal-110-Q2.csv};
\addplot [my line=brown!30!white] table [x=time, y=mechanics_iterations_linear_avg, col sep=comma] {data/beam/primal-011-Q2.csv};
\addplot [my line=green!30!white!90!black] table [x=time, y=mechanics_iterations_linear_avg, col sep=comma] {data/beam/primal-111-Q2.csv};
\legend {V, VE, EF, VEF}
\end{axis}
\end{tikzpicture}

%% file: tikz/primal/primal_swelling_Q1.tex
\newcommand{\indMin}{2}
\newcommand{\indMax}{2000}

\begin{tikzpicture}
\begin{axis}
[width=\textwidth,no markers, xlabel=Time $(s)$, ylabel=Average linear iterations, tick label style={font=\footnotesize}, legend style={fill opacity=0.6, font=\tiny,legend cell align=left, at={(0.98, 0.5)}, anchor=east}]
\addplot [my line=red!30!white] table [x=time, y=mechanics_iterations_linear_avg, col sep=comma] {data/swelling/primal-100-Q1.csv};
\addplot [my line=blue!30!white] table [x=time, y=mechanics_iterations_linear_avg, col sep=comma] {data/swelling/primal-110-Q1.csv};
\addplot [my line=brown!30!white] table [x=time, y=mechanics_iterations_linear_avg, col sep=comma] {data/swelling/primal-011-Q1.csv};
\addplot [my line=green!30!white!90!black] table [x=time, y=mechanics_iterations_linear_avg, col sep=comma] {data/swelling/primal-111-Q1.csv};
\legend {V, VE, EF, VEF}
\end{axis}
\end{tikzpicture}

%% file: tikz/primal/primal_swelling_Q2.tex
\begin{tikzpicture}
\begin{axis}
[width=\textwidth, no markers, xlabel=Time $(s)$, ylabel=Average linear iterations, tick label style={font=\footnotesize}, legend style={fill opacity=0.6, font=\tiny,legend cell align=left, at={(0.98, 0.5)}, anchor=east}]
\addplot [my line=red!30!white] table [x=time, y=mechanics_iterations_linear_avg, col sep=comma] {data/swelling/primal-100-Q2.csv};
\addplot [my line=blue!30!white] table [x=time, y=mechanics_iterations_linear_avg, col sep=comma] {data/swelling/primal-110-Q2.csv};
\addplot [my line=brown!30!white] table [x=time, y=mechanics_iterations_linear_avg, col sep=comma] {data/swelling/primal-011-Q2.csv};
\addplot [my line=green!30!white!90!black] table [x=time, y=mechanics_iterations_linear_avg, col sep=comma] {data/swelling/primal-111-Q2.csv};
\legend {V, VE, EF, VEF}
\end{axis}
\end{tikzpicture}

%% file: tikz/primal/primal_realistic_Q1.tex
\newcommand{\indMin}{2}
\newcommand{\indMax}{2000}

\begin{tikzpicture}
\begin{axis}
[width=\textwidth, no markers, xlabel=Time $(s)$, ylabel=Average linear iterations, tick label style={font=\footnotesize}, legend style={fill opacity=0.6, font=\tiny,legend cell align=left}, legend pos=north west]
\addplot [my line=red!30!white] table [x=time, y=mechanics_iterations_linear_avg, col sep=comma] {data/realistic/primal-100-Q1.csv};
\addplot [my line=blue!30!white] table [x=time, y=mechanics_iterations_linear_avg, col sep=comma] {data/realistic/primal-110-Q1.csv};
\addplot [my line=brown!30!white] table [x=time, y=mechanics_iterations_linear_avg, col sep=comma] {data/realistic/primal-011-Q1.csv};
\addplot [my line=green!30!white!90!black] table [x=time, y=mechanics_iterations_linear_avg, col sep=comma] {data/realistic/primal-111-Q1.csv};
\legend {V, VE, EF, VEF}
\end{axis}
\end{tikzpicture}

%% file: tikz/primal/primal_realistic_Q2.tex
\newcommand{\indMin}{2}
\newcommand{\indMax}{2000}

\begin{tikzpicture}
\begin{axis}
[width=\textwidth, no markers, xlabel=Time $(s)$, ylabel=Average linear iterations, tick label style={font=\footnotesize}, legend style={fill opacity=0.6, font=\tiny,legend cell align=left}, legend pos=north west]
\addplot [my line=red!30!white] table [x=time, y=mechanics_iterations_linear_avg, col sep=comma] {data/realistic/primal-100-Q2.csv};
\addplot [my line=blue!30!white] table [x=time, y=mechanics_iterations_linear_avg, col sep=comma] {data/realistic/primal-110-Q2.csv};
\addplot [my line=brown!30!white] table [x=time, y=mechanics_iterations_linear_avg, col sep=comma] {data/realistic/primal-011-Q2.csv};
\addplot [my line=green!30!white!90!black] table [x=time, y=mechanics_iterations_linear_avg, col sep=comma] {data/realistic/primal-111-Q2.csv};
\legend {V, VE, EF, VEF}
\end{axis}
\end{tikzpicture}

%% file: tikz/strong/strong_scaling_beam_time.tex
    \begin{tikzpicture}
    \begin{loglogaxis}
        [width=\textwidth,
        xlabel=CPUs,
        ylabel=Time (s),
         legend style={fill opacity=0.6, font=\tiny,legend cell align=left}, 
         legend columns=2
        ]
        
        \addplot [BDDCP1style] coordinates { (2,2829.713) (4,1451.548) (8,751.538) (16,504.714) (32,283.318) (64,186.046) (128,79.930) };
                                
        \addplot [AMGP1style] coordinates { (2,1965.827) (4,1055.379) (8,575.939) (16,328.499) (32,226.088) (64,141.423) (128,117.802)};

        \addplot [BDDCP2style] coordinates { (2,5357.591) (4,2739.102) (8,1332.132) (16,767.195) (32,517.437) (64,309.446) (128,140.030) };
                                
        \addplot [AMGP2style] coordinates { (2,4278.399) (4,2269.168) (8,1290.848) (16,714.628) (32,677.540) (64,513.515) (128,328.896) };

                                

                                

        \addplot [color=black] coordinates {
                                (2,  5000*1/2)
                                (4,  5000*1/4)
                                (8,  5000*1/8)
                                (16, 5000*1/16)
                                (32, 5000*1/32)
                                (64, 5000*1/64)
                                (128,5000*1/128)};
    \legend{\Qone BDDC, \Qone AMG, \Qtwo BDDC, \Qtwo AMG, Linear}
    \end{loglogaxis}
    \end{tikzpicture}
    









%% file: tikz/strong/strong_scaling_beam_iterations.tex
    \begin{tikzpicture}
    \begin{semilogxaxis}
        [width=\textwidth,
        xlabel=CPUs,
        ylabel=Total linear iterations,
         legend style={fill opacity=0.6, font=\tiny,legend cell align=left}, 
         legend columns=2, 
         legend pos=north west
        ]
        \addplot [BDDCP1style] coordinates {(2,41) (4,76) (8,175) (16,393) (32,488) (64,583) (128,488)};
                                
        \addplot [AMGP1style] coordinates {(2,346) (4,350) (8,349) (16,352) (32,361) (64,350) (128,365)};

        \addplot [BDDCP2style] coordinates {(2,51) (4,96) (8,226) (16,544) (32,697) (64,784) (128,671)};
                                
        \addplot [AMGP2style] coordinates {(2,624) (4,609) (8,631) (16,622) (32,640) (64,1071) (128,1261)};        
                                

                                

    \legend{\Qone BDDC, \Qone AMG, \Qtwo BDDC, \Qtwo AMG}
    \end{semilogxaxis}
    \end{tikzpicture}
    









%% file: tikz/strong/strong_scaling_swelling_time.tex
    \begin{tikzpicture}
    \begin{loglogaxis}
        [width=\textwidth,
        xlabel=CPUs,
        ylabel=Time (s),
        legend style={fill opacity=0.6, font=\tiny,legend cell align=left}, 
        legend columns=2
        ]

        \addplot [BDDCP1style] coordinates {(2,3015.024) (4,1575.033) (8,830.097) (16,490.530) (32,456.762) (64,224.238) (128,162.259)};
        \addplot [AMGP1style] coordinates {(2,2830.176) (4,1570.163) (8,896.732) (16,528.303) (32,384.262) (64,271.483) (128,190.100)};
                                
        \addplot [BDDCP2style] coordinates {(2,3864.443) (4,2040.645) (8,1078.262) (16,599.991) (32,429.986) (64,250.823) (128,172.612)};
        \addplot [AMGP2style] coordinates {(2,4130.198) (4,2260.834) (8,1313.164) (16,891.929) (32,745.571) (64,439.120) (128,275.780)};

                                

        \addplot [color=black] coordinates {
                                (2,  5000*1/2)
                                (4,  5000*1/4)
                                (8,  5000*1/8)
                                (16, 5000*1/16)
                                (32, 5000*1/32)
                                (64, 5000*1/64)
                                (128,5000*1/128)};
    \legend{\Qone BDDC, \Qone AMG, \Qtwo BDDC, \Qtwo AMG, Linear}
    \end{loglogaxis}
    \end{tikzpicture}
    









%% file: tikz/strong/strong_scaling_swelling_iterations.tex
    \begin{tikzpicture}
    \begin{semilogxaxis}
        [width=\textwidth,
        xlabel=CPUs,
        ylabel=Total linear iterations,
         legend style={fill opacity=0.6, font=\tiny,legend cell align=left, at={(0.98, 0.55)}, anchor=east}
        ]
        
        \addplot [BDDCP1style] coordinates {(2,431) (4,511) (8,540) (16,719) (32,798) (64,838) (128,766)};
        
        \addplot [AMGP1style] coordinates {(2,271) (4,271) (8,273) (16,270) (32,273) (64,273) (128,273)};
                                
        \addplot [BDDCP2style] coordinates {(2,565) (4,690) (8,687) (16,749) (32,780) (64,724) (128,691)};
        
        \addplot [AMGP2style] coordinates {(2,445) (4,444) (8,443) (16,443) (32,443) (64,442) (128,443)};
        
                                

    \legend{\Qone BDDC, \Qone AMG, \Qtwo BDDC, \Qtwo AMG}
    \end{semilogxaxis}
    \end{tikzpicture}
    

%% file: tikz/strong/strong_scaling_contraction_time.tex
    \begin{tikzpicture}
    \begin{loglogaxis}
        [width=\textwidth,
        xlabel=CPUs,
        ylabel=Time (s),
        legend style={fill opacity=0.6, font=\tiny,legend cell align=left},
        legend columns=2 
        ]
        
        \addplot [BDDCP1style] coordinates {(2,3502.742) (4,1915.705) (8,1008.466) (16,600.301) (32,483.097) (64,281.492) (128,150.714)};  

        \addplot [AMGP1style] coordinates {(2,3392.544) (4,1854.422) (8,1046.480) (16,630.686) (32,507.715) (64,303.874) (128,203.381)};  
                            
        \addplot [BDDCP2style] coordinates {(2,4209.062) (4,2267.183) (8,1185.835) (16,706.704) (32,550.101) (64,297.043) (128,186.413)};  
        \addplot [AMGP2style] coordinates {(2,4602.675) (4,2493.259) (8,1393.627) (16,828.889) (32,772.616) (64,434.825) (128,285.030)}; 

        \addplot [color=black] coordinates {
                                (2,  5000*1/2)
                                (4,  5000*1/4)
                                (8,  5000*1/8)
                                (16, 5000*1/16)
                                (32, 5000*1/32)
                                (64, 5000*1/64)
                                (128,5000*1/128)};
    \legend{\Qone BDDC, \Qone AMG, \Qtwo BDDC, \Qtwo AMG, Linear}
    \end{loglogaxis}
    \end{tikzpicture}
    









%% file: tikz/strong/strong_scaling_contraction_iterations.tex
    \begin{tikzpicture}
    \begin{semilogxaxis}
        [width=\textwidth,
        xlabel=CPUs,
        ylabel=Total linear iterations,
        legend style={fill opacity=0.6, 
                        font=\tiny,
                        legend cell align=left,
                        at={(0.98,0.55)}, anchor=east}
        ]

        \addplot [BDDCP1style] coordinates {(2,389) (4,512) (8,540) (16,693) (32,724) (64,774) (128,760)};  

        \addplot [AMGP1style] coordinates {(2,261) (4,262) (8,261) (16,261) (32,265) (64,265) (128,264)};  
                                
        \addplot [BDDCP2style] coordinates {(2,505) (4,664) (8,666) (16,705) (32,722) (64,682) (128,664)};  
        \addplot [AMGP2style] coordinates {(2,437) (4,436) (8,435) (16,431) (32,429) (64,426) (128,427)};
        


    \legend{\Qone BDDC, \Qone AMG, \Qtwo BDDC, \Qtwo AMG}
    \end{semilogxaxis}
    \end{tikzpicture}
    

%% file: tikz/stronger/stronger_scaling_beam_time_Q1.tex
    \begin{tikzpicture}
    \begin{loglogaxis}
        [width=\textwidth,
        ylabel=Time (s),
        legend style={fill opacity=0.6, font=\tiny,legend cell align=left}, 
        legend columns=2
        ]
        \addplot [BDDCP1style] coordinates {
                                (32,1505.482)
                                (64,1001.786)
                                (128,407.101)
                                (256,133.550)
                                (512,76.062)
                                (1024,104.358)};
        \addplot [AMGP1style] coordinates {
                                (32,709.552)
                                (64,430.385)
                                (128,413.378)
                                (256,301.695)
                                (512,118.174)
                                (1024,90.469)};
        \addplot [BDDCMLP1style] coordinates {
                                (32,1456.862)
                                (64,666.267)
                                (128,365.919)
                                (256,212.330)
                                (512,99.666)
                                (1024,65.315)
                                };
        \addplot [color=black] coordinates {
                                (32, 30000*1/32)
                                (64, 30000*1/64)
                                (128,30000*1/128)
                                (256,30000*1/256)
                                (512,30000*1/512)};
    \legend{\Qone BDDC, \Qone AMG, \Qone BDDC-ML, Linear}
    \end{loglogaxis}
    \end{tikzpicture}
    









%% file: tikz/stronger/stronger_scaling_beam_iterations_Q1.tex
    \begin{tikzpicture}
    \begin{semilogxaxis}
        [width=\textwidth,
        ylabel=Average linear iterations,
        legend style={fill opacity=0.6, font=\tiny,legend cell align=left}, 
        legend pos=north west
        ]
        
        \addplot [BDDCP1style] coordinates {
                               (32,228) (64,441) (128,375) (256,238) (512,232) (1024,242)};
                                
        \addplot [AMGP1style] coordinates {
                                (32,137) (64,156) (128,299) (256,155) (512,140) (1024,142)};

        \addplot [BDDCMLP1style] coordinates {
                                (32,220) (64,303) (128,366) (256,602) (512,603) (1024,890)};

    \legend{\Qone BDDC, \Qone AMG, \Qone BDDC-ML}
    \end{semilogxaxis}
    \end{tikzpicture}
    









%% file: tikz/stronger/stronger_scaling_beam_time_Q2.tex
    \begin{tikzpicture}
    \begin{loglogaxis}
        [width=\textwidth,
        ylabel=Time (s),
         legend style={fill opacity=0.6, font=\tiny,legend cell align=left},
         legend columns=2
        ]

        \addplot [BDDCP2style] coordinates {
                                (32,2344.275)
                                (64,1876.059)
                                (128,797.997)
                                (256,254.602)
                                (512,152.210)
                                (1024,164.870)};
                                
        \addplot [AMGP2style] coordinates { (32,1811.482) (64,1881.091) (128,1451.428) (256,973.231) (512,553.363) (1024,531.543)
                                };
                                
        \addplot [BDDCMLP2style] coordinates {
                                (32,2532.437)
                                (64,1192.754)
                                (128,627.026)
                                (256,360.314)
                                (512,174.769)
                                (1024,105.418)
                                };

        \addplot [color=black] coordinates {
                                (32, 50000*1/32)
                                (64, 50000*1/64)
                                (128,50000*1/128)
                                (256,50000*1/256)
                                (512,50000*1/512)};
    \legend{\Qtwo BDDC, \Qtwo AMG, \Qtwo BDDC-ML, Linear}
    \end{loglogaxis}
    \end{tikzpicture}
    

%% file: tikz/stronger/stronger_scaling_beam_iterations_Q2.tex
    \begin{tikzpicture}
    \begin{semilogxaxis}
        [width=\textwidth,
        ylabel=Average linear iterations,
        legend style={fill opacity=0.6, font=\tiny,legend cell align=left},
        legend pos=north west
        ]
        
        \addplot [BDDCP2style] coordinates { (32,324) (64,627) (128,557) (256,362) (512,388) (1024,386) };
                                
        \addplot [AMGP2style] coordinates { (32,278) (64,746) (128,1034) (256,1099) (512,1622) (1024,1691) };
        
        \addplot [BDDCMLP2style] coordinates { (32,313) (64,399) (128,501) (256,816) (512,840) (1024,1211) };

    \legend{\Qtwo BDDC, \Qtwo AMG, \Qtwo BDDC-ML}
    \end{semilogxaxis}
    \end{tikzpicture}
  

%% file: tikz/stronger/stronger_scaling_swelling_time_Q1.tex
    \begin{tikzpicture}
    \begin{loglogaxis}
        [width=\textwidth,
        ylabel=Time (s),
         legend style={fill opacity=0.6, font=\tiny,legend cell align=left}
        ]

        \addplot [BDDCP1style] coordinates { (32,1909.709) (64,897.746) (128,493.733) (256,245.393) (512,236.053) (1024,538.329) };  
        
        \addplot [AMGP1style] coordinates {
                                (32,1466.151)
                                (64,959.444)
                                (128,788.265)
                                (256,677.673)
                                (512,206.759)
                                (1024,150.724)
                                }; 
                                
        \addplot [BDDCMLP1style] coordinates {
                                (32,1916.671)
                                (64,1041.185)
                                (128,584.721)
                                (256,431.377)
                                (512,228.296)
                                (1024,226.272)
                                }; 

        \addplot [color=black] coordinates {
                                (32, 50000*1/32)
                                (64, 50000*1/64)
                                (128,50000*1/128)
                                (256,50000*1/256)
                                (512,50000*1/512)};
        \legend{\Qone BDDC, \Qone AMG, \Qone BDDC-ML, Linear}
    \end{loglogaxis}
    \end{tikzpicture}
    









%% file: tikz/stronger/stronger_scaling_swelling_iterations_Q1.tex
    \begin{tikzpicture}
    \begin{loglogaxis}
        [width=\textwidth,
        ylabel=Average linear iterations,
         legend style={fill opacity=0.6, font=\tiny,legend cell align=left}, 
         legend pos=north west
        ]

        \addplot [BDDCP1style] coordinates {
                               (32,524) (64,588) (128,546) (256,541) (512,490) (1024,444)};
                                
        \addplot [AMGP1style] coordinates { (32,174) (64,173) (128,174) (256,174) (512,175) (1024,178)};

        \addplot [BDDCMLP1style] coordinates { (32,522) (64,716) (128,774) (256,1472) (512,1427) (1024,3252)};

    \legend{\Qone BDDC, \Qone AMG, \Qone BDDC-ML}
    \end{loglogaxis}
    \end{tikzpicture}
    

%% file: tikz/stronger/stronger_scaling_swelling_time_Q2.tex
    \begin{tikzpicture}
    \begin{loglogaxis}
        [width=\textwidth,
        ylabel=Time (s),
         legend style={fill opacity=0.6, font=\tiny,legend cell align=left}
        ]

        \addplot [BDDCP2style] coordinates {
                                (32,2516.032)
                                (64,1159.916)
                                (128,623.832)
                                (256,345.394)
                                (512,316.504)
                                (1024,691.804)};  
        \addplot [AMGP2style] coordinates {
                                (32,2291.441)
                                (64,1429.323)
                                (128,1023.769)
                                (256,856.772)
                                (512,266.129)
                                (1024,171.158)
                                };  
                                
        \addplot [BDDCMLP2style] coordinates { (32,2605.377) (64,1435.504) (128,763.625) (256,459.409) (512,305.853) (1024,250.795) };

        \addplot [color=black] coordinates {
                                (32, 50000*1/32)
                                (64, 50000*1/64)
                                (128,50000*1/128)
                                (256,50000*1/256)
                                (512,50000*1/512)};
        \legend{\Qtwo BDDC, \Qtwo AMG, \Qtwo BDDC-ML, Linear}
    \end{loglogaxis}
    \end{tikzpicture}

%% file: tikz/stronger/stronger_scaling_swelling_iterations_Q2.tex
    \begin{tikzpicture}
    \begin{loglogaxis}
        [width=\textwidth,
        ylabel=Average linear iterations,
         legend style={fill opacity=0.6, font=\tiny,legend cell align=left},
         legend pos=north west
        ]

        \addplot [BDDCP2style] coordinates { (32,324) (64,627) (128,557) (256,362) (512,388) (1024,386) };
                                
        \addplot [AMGP2style] coordinates { (32,279) (64,280) (128,280) (256,279) (512,275) (1024,280) };
        
        \addplot [BDDCMLP2style] coordinates { (32,628) (64,910) (128,921) (256,1251) (512,1731) (1024,2685) };

    \legend{\Qtwo BDDC, \Qtwo AMG, \Qtwo BDDC-ML}
    \end{loglogaxis}
    \end{tikzpicture}
    

%% file: tikz/stronger/stronger_scaling_contraction_time_Q1.tex
    \begin{tikzpicture}
    \begin{loglogaxis}
        [width=\textwidth,
        ylabel=Time (s),
         legend style={fill opacity=0.6, font=\tiny,legend cell align=left}
        ]
        
        \addplot [BDDCP1style] coordinates {
                                (32,1931.875)
                                (64,951.858)
                                (128,485.244)
                                (256,270.793)
                                (512,251.502)
                                (1024,598.210)};  

        \addplot [AMGP1style] coordinates {
                                (32,1772.135)
                                (64,1079.338)
                                (128,872.050)
                                (256,776.572)
                                (512,235.105)
                                (1024,164.619)
                                };           

        \addplot [BDDCMLP1style] coordinates {
                                (32,2066.918)
                                (64,1050.427)
                                (128,580.401)
                                (256,425.828)
                                (512,207.855)
                                (1024,138.964)
                                };

        \addplot [color=black] coordinates {
                                (32, 30000*1/32)
                                (64, 30000*1/64)
                                (128,30000*1/128)
                                (256,30000*1/256)
                                (512,30000*1/512)};
    \legend{\Qone BDDC, \Qone AMG, \Qone BDDC-ML, Linear}
    \end{loglogaxis}
    \end{tikzpicture}
    









%% file: tikz/stronger/stronger_scaling_contraction_iterations_Q1.tex
    \begin{tikzpicture}
    \begin{semilogxaxis}
        [width=\textwidth,
        ylabel=Average linear iterations,
         legend style={fill opacity=0.6, font=\tiny,legend cell align=left}, 
         legend pos=north west
        ]
        
        \addplot [BDDCP1style] coordinates {(32,464) (64,520) (128,504) (256,519) (512,479) (1024,457)};
                                
        \addplot [AMGP1style] coordinates { (32,181) (64,181) (128,178) (256,180) (512,197) (1024,181)};

        \addplot [BDDCMLP1style] coordinates {
                                (32,468) (64,602) (128,648) (256,1276) (512,1153) (1024,1729)};

    \legend{\Qone BDDC, \Qone AMG, \Qone BDDC-ML}
    \end{semilogxaxis}
    \end{tikzpicture}
    

%% file: tikz/stronger/stronger_scaling_contraction_time_Q2.tex
    \begin{tikzpicture}
    \begin{loglogaxis}
        [width=\textwidth,
        ylabel=Time (s),
         legend style={fill opacity=0.6, font=\tiny,legend cell align=left}
        ]

        \addplot [BDDCP2style] coordinates {
                                (32,2634.190)
                                (64,1202.147)
                                (128,636.893)
                                (256,373.173)
                                (512,329.862)
                                (1024,759.256)};  
        \addplot [AMGP2style] coordinates {
                                (32,2522.257)
                                (64,1551.130)
                                (128,1093.390)
                                (256,934.331)
                                (512,295.707)
                                (1024,191.135)};  

        \addplot [BDDCMLP2style] coordinates {
                                (32,2686.669)
                                (64,1275.709)
                                (128,739.828)
                                (256,400.107)
                                (512,230.980)
                                (1024,159.754)};  

        \addplot [color=black] coordinates {
                                (32, 30000*1/32)
                                (64, 30000*1/64)
                                (128,30000*1/128)
                                (256,30000*1/256)
                                (512,30000*1/512)};
    \legend{\Qtwo BDDC, \Qtwo AMG, \Qtwo BDDC-ML, Linear}
    \end{loglogaxis}
    \end{tikzpicture}

%% file: tikz/stronger/stronger_scaling_contraction_iterations_Q2.tex
    \begin{tikzpicture}
    \begin{semilogxaxis}
        [width=\textwidth,
        ylabel=Average linear iterations,
         legend style={fill opacity=0.6, font=\tiny,legend cell align=left},
         legend pos=north west
        ]

        \addplot [BDDCP2style] coordinates { (32,563) (64,564) (128,584) (256,641) (512,615) (1024,609) };
                                
        \addplot [AMGP2style] coordinates { (32,287) (64,288) (128,287) (256,287) (512,287) (1024,320) };
        
        \addplot [BDDCMLP2style] coordinates { (32,559) (64,619) (128,735) (256,927) (512,1127) (1024,1668) };

    \legend{\Qtwo BDDC, \Qtwo AMG, \Qtwo BDDC-ML}
    \end{semilogxaxis}
    \end{tikzpicture}
    

%% file: tikz/EM/small-p1.tex
    \newcommand{\tagcolor}{yellow!20!white}
    \newcommand{\tzero}{0}
    \newcommand{\ttagone}{0.075}
    \newcommand{\tone}{0.13}
    \newcommand{\ttagtwo}{0.25}
    \newcommand{\ttwo}{0.31}
    \newcommand{\ttagthree}{0.44}
    \newcommand{\tthree}{0.475}
    \newcommand{\ttagfour}{0.57}
    \newcommand{\tfour}{0.80}
    \newcommand{\ytag}{120}
    \newcommand{\ytop}{200}
    \newcommand{\timetag}{15}
    \newcommand{\timetop}{20}
    \centering
    \begin{subfigure}{\textwidth}
    \begin{semilogyplot}[height=7cm, width=\textwidth, legend pos=north east,  legend style={at={(0.9, 0.9)}, anchor=north east}]{Time}{Average linear its.}
        \addplot+[AMGP1style, mark=none] table [x=time, y=mechanics_iterations_linear_avg, col sep=comma, select coords between index={2}{160}] {data/EM/electromechanics-sr-small.csv};
        \addplot+[BDDCP1style, mark=none] table [x=time, y=mechanics_iterations_linear_avg, col sep=comma, select coords between index={2}{160}] {data/EM/electromechanics-sr-bddc-small.csv};
        \addplot+[fill, mark=none, opacity=0.1, draw=none] coordinates 
		{(\tzero,0) (\tone,0) (\tone,\ytop) (0,\ytop)} \closedcycle;
		\node[fill=\tagcolor, text centered] 
		at (\ttagone,\ytag) {\small IC};
		\addplot+[fill, mark=none, opacity=0.1, draw=none] coordinates 
		{(\tone,0) (\ttwo,0) (\ttwo,\ytop) (\tone,\ytop)} \closedcycle;
		\node[fill=\tagcolor, text centered] 
		at (\ttagtwo,\ytag) {\small C};
		\addplot+[fill, mark=none, opacity=0.1, draw=none] coordinates 
		{(\ttwo,0) (\tthree,0) (\tthree,\ytop) (\ttwo,\ytop)} \closedcycle;
		\node[fill=\tagcolor, text centered] 
		at (\ttagthree,\ytag) {\small IR};
		\addplot+[fill, mark=none, opacity=0.1, draw=none] coordinates 
		{(\tthree,0) (\tfour,0) (\tfour,\ytop) (\tthree,\ytop)} \closedcycle;
		\node[fill=\tagcolor, text centered] 
		at (\ttagfour,\ytag) {\small R};
        \legend{AMG, BDDC}
    \end{semilogyplot}
    \end{subfigure}
    \begin{subfigure}{\textwidth}
    \begin{semilogyplot}[height=7cm,  legend style={at={(0.9, 0.9)}, anchor=north east}]{Time}{Solution time}
        \addplot+[AMGP1style, mark=none] table [x=time, y=mechanics_solving_time, col sep=comma, select coords between index={2}{160}] {data/EM/electromechanics-sr-small.csv};
        \addplot+[BDDCP1style, mark=none] table [x=time, y=mechanics_solving_time, col sep=comma, select coords between index={2}{160}] {data/EM/electromechanics-sr-bddc-small.csv};
        \addplot+[fill, mark=none, opacity=0.1, draw=none] coordinates 
		{(\tzero,0) (\tone,0) (\tone,\timetop) (0,\timetop)} \closedcycle;
		\node[fill=\tagcolor, text centered] 
		at (\ttagone,\timetag) {\small IC};
		\addplot+[fill, mark=none, opacity=0.1, draw=none] coordinates 
		{(\tone,0) (\ttwo,0) (\ttwo,\timetop) (\tone,\timetop)} \closedcycle;
		\node[fill=\tagcolor, text centered] 
		at (\ttagtwo,\timetag) {\small C};
		\addplot+[fill, mark=none, opacity=0.1, draw=none] coordinates 
		{(\ttwo,0) (\tthree,0) (\tthree,\timetop) (\ttwo,\timetop)} \closedcycle;
		\node[fill=\tagcolor, text centered] 
		at (\ttagthree,\timetag) {\small IR};
		\addplot+[fill, mark=none, opacity=0.1, draw=none] coordinates 
		{(\tthree,0) (\tfour,0) (\tfour,\timetop) (\tthree,\timetop)} \closedcycle;
		\node[fill=\tagcolor, text centered] 
		at (\ttagfour,\timetag) {\small R};
        \legend{AMG, BDDC}
    \end{semilogyplot}
    \end{subfigure}

%% file: tikz/EM/large-p1.tex
    \newcommand{\tagcolor}{yellow!20!white}
    \newcommand{\tzero}{0}
    \newcommand{\ttagone}{0.075}
    \newcommand{\tone}{0.13}
    \newcommand{\ttagtwo}{0.23}
    \newcommand{\ttwo}{0.31}
    \newcommand{\ttagthree}{0.42}
    \newcommand{\tthree}{0.475}
    \newcommand{\ttagfour}{0.57}
    \newcommand{\tfour}{0.80}
    \newcommand{\ytag}{300}
    \newcommand{\ytop}{750}
    \newcommand{\timetag}{90}
    \newcommand{\timetop}{220}
    \centering
    \begin{subfigure}{0.48\textwidth}
    \begin{semilogyplot}[height=7cm, legend pos=north east,  legend style={at={(0.9, 0.5)}, anchor=east}]{Time}{Average linear its.}
        \addplot+[AMGP1style, mark=none] table [x=time, y=mechanics_iterations_linear_avg, col sep=comma, select coords between index={2}{62}] {data/EM/electromechanics-sr.csv};
        \addplot+[BDDCP1style, mark=none] table [x=time, y=mechanics_iterations_linear_avg, col sep=comma, select coords between index={2}{62}] {data/EM/electromechanics-sr-bddc.csv};
        \addplot+[fill, mark=none, opacity=0.1, draw=none] coordinates 
		{(\tzero,0) (\tone,0) (\tone,\ytop) (0,\ytop)} \closedcycle;
		\node[fill=\tagcolor, text centered] 
		at (\ttagone,\ytag) {\small IC};
		\addplot+[fill, mark=none, opacity=0.1, draw=none] coordinates 
		{(\tone,0) (\ttwo,0) (\ttwo,\ytop) (\tone,\ytop)} \closedcycle;
		\node[fill=\tagcolor, text centered] 
		at (\ttagtwo,\ytag) {\small C};
        \legend{AMG, BDDC}
    \end{semilogyplot}
    \end{subfigure}
    \begin{subfigure}{0.48\textwidth}
    \begin{semilogyplot}[height=7cm,  legend style={at={(0.9, 0.5)}, anchor=north east}]{Time}{Solution time}
        \addplot+[AMGP1style, mark=none] table [x=time, y=mechanics_solving_time, col sep=comma, select coords between index={2}{62}] {data/EM/electromechanics-sr.csv};
        \addplot+[BDDCP1style, mark=none] table [x=time, y=mechanics_solving_time, col sep=comma, select coords between index={2}{62}] {data/EM/electromechanics-sr-bddc.csv};
        \addplot+[fill, mark=none, opacity=0.1, draw=none] coordinates 
		{(\tzero,0) (\tone,0) (\tone,\timetop) (0,\timetop)} \closedcycle;
		\node[fill=\tagcolor, text centered] 
		at (\ttagone,\timetag) {\small IC};
		\addplot+[fill, mark=none, opacity=0.1, draw=none] coordinates 
		{(\tone,0) (\ttwo,0) (\ttwo,\timetop) (\tone,\timetop)} \closedcycle;
		\node[fill=\tagcolor, text centered] 
		at (\ttagtwo,\timetag) {\small C};
        \legend{AMG, BDDC}
    \end{semilogyplot}
    \end{subfigure}

%% file: main.bbl
\newcommand{\etalchar}[1]{$^{#1}$}
\begin{thebibliography}{FPMM{\etalchar{+}}19}

\bibitem[AANQ11]{ambrosi2011electromechanical}
D.~Ambrosi, G.~Arioli, F.~Nobile, and A.M. Quarteroni.
\newblock Electromechanical coupling in cardiac dynamics: the active strain
  approach.
\newblock {\em SIAM Journal on Applied Mathematics}, 71(2):605--621, 2011.

\bibitem[Ada02]{adams2002evaluation}
M.~Adams.
\newblock Evaluation of three unstructured multigrid methods on 3d finite
  element problems in solid mechanics.
\newblock {\em International Journal for Numerical Methods in Engineering},
  55(5):519--534, 2002.

\bibitem[ADLK00]{amestoy2000mumps}
P.R. Amestoy, I.S. Duff, J.-Y. L’Excellent, and J.~Koster.
\newblock Mumps: a general purpose distributed memory sparse solver.
\newblock In {\em International Workshop on Applied Parallel Computing}, pages
  121--130. Springer, 2000.

\bibitem[ANL{\etalchar{+}}16]{AUGUSTIN2016622}
C.M. Augustin, A.~Neic, M.~Liebmann, A.J. Prassl, S.A. Niederer, G.~Haase, and
  G.~Plank.
\newblock Anatomically accurate high resolution modeling of human whole heart
  electromechanics: A strongly scalable algebraic multigrid solver method for
  nonlinear deformation.
\newblock {\em Journal of Computational Physics}, 305:622--646, 2016.

\bibitem[AP95]{ambrosetti1995primer}
A.~Ambrosetti and G.~Prodi.
\newblock {\em A primer of nonlinear analysis}.
\newblock Number~34. Cambridge University Press, 1995.

\bibitem[APF{\etalchar{+}}22]{lifex}
P.C. Africa, R.~Piersanti, M.~Fedele, L.~Dede, and A.M. Quarteroni.
\newblock lifex--heart module: a high-performance simulator for the cardiac
  function package 1: Fiber generation.
\newblock {\em arXiv preprint arXiv:2201.03303}, 2022.

\bibitem[BAA{\etalchar{+}}21]{petsc-user-ref}
S.~Balay, S.~Abhyankar, M.F. Adams, J.~Brown, P.~Brune, K.~Buschelman,
  L.~Dalcin, A.~Dener, V.~Eijkhout, W.~Gropp, D.~Karpeyev, D.~Kaushik,
  M.~Knepley, D.~May, L.~Curfman~McInnes, R.~Mills, T.~Munson, K.~Rupp,
  P.~Sanan, B.~Smith, S.~Zampini, H.~Zhang, and H.~Zhang.
\newblock {PETS}c users manual.
\newblock Technical Report ANL-95/11 - Revision 3.13, Argonne National
  Laboratory, 2021.

\bibitem[BBPT12]{Bayer20122243}
J.D. Bayer, R.C. Blake, G.~Plank, and N.A. Trayanova.
\newblock A novel rule-based algorithm for assigning myocardial fiber
  orientation to computational heart models.
\newblock {\em Annals of Biomedical Engineering}, 40(10):2243--2254, 2012.

\bibitem[BHK07]{dealii}
W.~Bangerth, R.~Hartmann, and G.~Kanschat.
\newblock deal. ii—a general-purpose object-oriented finite element library.
\newblock {\em ACM Transactions on Mathematical Software (TOMS)}, 33(4):24--es,
  2007.

\bibitem[BTB06]{brezina2006parallel}
M.~Brezina, C.~Tong, and R.~Becker.
\newblock Parallel algebraic multigrids for structural mechanics.
\newblock {\em SIAM Journal on Scientific Computing}, 27(5):1534--1554, 2006.

\bibitem[CD04]{carstensen2004priori}
C.~Carstensen and G.~Dolzmann.
\newblock An a priori error estimate for finite element discretizations in
  nonlinear elasticity for polyconvex materials under small loads.
\newblock {\em Numerische Mathematik}, 97(1):67--80, 2004.

\bibitem[CFPS14]{franzone2014mathematical}
P.~Colli~Franzone, L.F. Pavarino, and S.~Scacchi.
\newblock {\em Mathematical Cardiac Electrophysiology}, volume~13.
\newblock Springer, 2014.

\bibitem[CFPS18]{colli2018numerical}
P.~Colli~Franzone, L.F. Pavarino, and S.~Scacchi.
\newblock A numerical study of scalable cardiac electro-mechanical solvers on
  {HPC} architectures.
\newblock {\em Frontiers in Physiology}, 9:268, 2018.

\bibitem[CP08]{chevalier2008pt}
C.~Chevalier and F.~Pellegrini.
\newblock {PT}-{S}cotch: A tool for efficient parallel graph ordering.
\newblock {\em Parallel computing}, 34(6-8):318--331, 2008.

\bibitem[Dav07]{davis2007umfpack}
T.A. Davis.
\newblock Umfpack version 5.2. 0 user guide.
\newblock {\em University of Florida}, 25, 2007.

\bibitem[DN10]{davis2010algorithm}
T.A. Davis and E.P. Natarajan.
\newblock Algorithm 907: Klu, a direct sparse solver for circuit simulation
  problems.
\newblock {\em ACM Transactions on Mathematical Software (TOMS)}, 37(3):1--17,
  2010.

\bibitem[Doh03]{dohrmann2003preconditioner}
C.R. Dohrmann.
\newblock A preconditioner for substructuring based on constrained energy
  minimization.
\newblock {\em SIAM Journal on Scientific Computing}, 25(1):246--258, 2003.

\bibitem[EMFTF10]{el2010iterative}
A.~El~Maliki, M.~Fortin, N.~Tardieu, and A.~Fortin.
\newblock Iterative solvers for 3d linear and nonlinear elasticity problems:
  Displacement and mixed formulations.
\newblock {\em International journal for numerical methods in engineering},
  83(13):1780--1802, 2010.

\bibitem[FLL{\etalchar{+}}01]{farhat2001feti}
C.~Farhat, M.~Lesoinne, P.~LeTallec, K.~Pierson, and D.~Rixen.
\newblock {FETI}-{DP}: a dual--primal unified {FETI} method—part i: A faster
  alternative to the two-level feti method.
\newblock {\em International journal for numerical methods in engineering},
  50(7):1523--1544, 2001.

\bibitem[FPMM{\etalchar{+}}19]{franceschini2019robust}
A.~Franceschini, V.A. Paduletto~Magri, G.~Mazzucco, N.~Spiezia, and C.~Janna.
\newblock A robust adaptive algebraic multigrid linear solver for structural
  mechanics.
\newblock {\em Computer Methods in Applied Mechanics and Engineering},
  352:389--416, 2019.

\bibitem[FY02]{falgout2002hypre}
R.D. Falgout and U.M. Yang.
\newblock hypre: A library of high performance preconditioners.
\newblock In {\em International Conference on Computational Science}, pages
  632--641. Springer, 2002.

\bibitem[GMW91]{Guccione1991}
J.M. Guccione, A.D. McCulloch, and L.K. Waldman.
\newblock Passive material properties of intact ventricular myocardium
  determined from a cylindrical model.
\newblock {\em Journal of biomechanical engineering}, 113(1):42--55, 1991.

\bibitem[GOS03]{griebel2003algebraic}
M.~Griebel, D.~Oeltz, and M.A. Schweitzer.
\newblock An algebraic multigrid method for linear elasticity.
\newblock {\em SIAM Journal on Scientific Computing}, 25(2):385--407, 2003.

\bibitem[Hol02]{Holzapfel2002489}
G.A. Holzapfel.
\newblock Nonlinear solid mechanics: A continuum approach for engineering
  science.
\newblock {\em Meccanica}, 37(4):489--490, 2002.

\bibitem[HW12]{heroux2012new}
M.A. Heroux and J.M. Willenbring.
\newblock A new overview of the trilinos project.
\newblock {\em Scientific Programming}, 20(2):83--88, 2012.

\bibitem[JCC20]{Jiang2020}
Y.~Jiang, R.~Chen, and X.-C. Cai.
\newblock A highly parallel implicit domain decomposition method for the
  simulation of the left ventricle on unstructured meshes.
\newblock {\em Computational Mechanics}, 66(6):1461--1475, 2020.

\bibitem[KGH{\etalchar{+}}22]{karabelas2022accurate}
E.~Karabelas, M.A.F. Gsell, G.~Haase, G.~Plank, and C.M. Augustin.
\newblock An accurate, robust, and efficient finite element framework with
  applications to anisotropic, nearly and fully incompressible elasticity.
\newblock {\em Computer Methods in Applied Mechanics and Engineering},
  394:114887, 2022.

\bibitem[KSK97]{karypis1997parmetis}
G.~Karypis, K.~Schloegel, and V.~Kumar.
\newblock Parmetis: Parallel graph partitioning and sparse matrix ordering
  library.
\newblock 1997.

\bibitem[Lev13]{levick2013introduction}
J.R. Levick.
\newblock {\em An introduction to cardiovascular physiology}.
\newblock Butterworth-Heinemann, 2013.

\bibitem[LGA{\etalchar{+}}15]{wrro93701}
S.~Land, V.~Gurev, S.~Arens, C.M. Augustin, L.~Baron, R.~Blake, C.~Bradley,
  S.~Castro, A.~Crozier, M.~Favino, T.E. Fastl, T.~Fritz, H.~Gao, A.~Gizzi,
  B.E. Griffith, D.E. Hurtado, R.~Krause, X.~Luo, M.P. Nash, S.~Pezzuto,
  G.~Plank, S.~Rossi, D.~Ruprecht, G.~Seemann, N.P. Smith, J.~Sundnes, J.J.
  Rice, N.~Trayanova, D.~Wang, Z.J Wang, and S.A. Niederer.
\newblock Verification of cardiac mechanics software: benchmark problems and
  solutions for testing active and passive material behaviour.
\newblock {\em Proceedings of the Royal Society A: Mathematical, Physical and
  Engineering Sciences}, 471(2184), December 2015.
\newblock {\copyright} 2015 The Authors. Published by the Royal Society under
  the terms of the Creative Commons Attribution License
  http://creativecommons.org/licenses/by/4.0/, which permits unrestricted use,
  provided the original author and source are credited.

\bibitem[Li05]{li2005overview}
X.S. Li.
\newblock An overview of superlu: Algorithms, implementation, and user
  interface.
\newblock {\em ACM Transactions on Mathematical Software (TOMS)},
  31(3):302--325, 2005.

\bibitem[LW06]{li2006feti}
J.~Li and O.B. Widlund.
\newblock {FETI}-{DP}, {BDDC}, and block {C}holesky methods.
\newblock {\em International journal for numerical methods in engineering},
  66(2):250--271, 2006.

\bibitem[MSD08]{mandel2008multispace}
J.~Mandel, B.~Soused{\'\i}k, and C.R. Dohrmann.
\newblock Multispace and multilevel {BDDC}.
\newblock {\em Computing}, 83(2-3):55--85, 2008.

\bibitem[MS{\v{S}}12]{mandel2012adaptive}
J.~Mandel, B.~Soused{\'\i}k, and J.~{\v{S}}{\'\i}stek.
\newblock Adaptive bddc in three dimensions.
\newblock {\em Mathematics and Computers in Simulation}, 82(10):1812--1831,
  2012.

\bibitem[PD17]{pechstein2017unified}
C.~Pechstein and C.R. Dohrmann.
\newblock A unified framework for adaptive bddc.
\newblock {\em Electron. Trans. Numer. Anal}, 46(273-336):3, 2017.

\bibitem[PHW{\etalchar{+}}19]{pfaller2019}
M.R. Pfaller, J.M. H{\"o}rmann, M.~Weigl, A.~Nagler, R.~Chabiniok,
  C.~Bertoglio, and W.A. Wall.
\newblock The importance of the pericardium for cardiac biomechanics: from
  physiology to computational modeling.
\newblock {\em Biomechanics and modeling in mechanobiology}, 18(2):503--529,
  2019.

\bibitem[PSZ15]{pavarinoSZ2015}
L.F. Pavarino, S.~Scacchi, and S.~Zampini.
\newblock Newton{\textendash}{K}rylov-{BDDC} solvers for nonlinear cardiac
  mechanics.
\newblock {\em Computer Methods in Applied Mechanics and Engineering},
  295:562--580, 2015.

\bibitem[QLRRB17]{quarteroni2017integrated}
A.M. Quarteroni, T.~Lassila, S.~Rossi, and R.~Ruiz-Baier.
\newblock Integrated heart—coupling multiscale and multiphysics models for
  the simulation of the cardiac function.
\newblock {\em Computer Methods in Applied Mechanics and Engineering},
  314:345--407, 2017.

\bibitem[QV08]{QuarteroniValli2008}
A.M. Quarteroni and A.~Valli.
\newblock {\em Numerical approximation of partial differential equations},
  volume~23.
\newblock Springer Science \& Business Media, 2008.

\bibitem[RSA{\etalchar{+}}20]{regazzoni2020cardiac}
F.~Regazzoni, M.~Salvador, P.C. Africa, M.~Fedele, L.~Dede, and A.M.
  Quarteroni.
\newblock A cardiac electromechanics model coupled with a lumped parameters
  model for closed-loop blood circulation. part {I}: model derivation.
\newblock {\em arXiv e-prints}, 2020.

\bibitem[Saa03]{saad2003iterative}
Y.~Saad.
\newblock {\em Iterative methods for sparse linear systems}.
\newblock SIAM, 2003.

\bibitem[Smi92]{smith1992optimal}
B.F. Smith.
\newblock An optimal domain decomposition preconditioner for the finite element
  solution of linear elasticity problems.
\newblock {\em SIAM Journal on Scientific and Statistical Computing},
  13(1):364--378, 1992.

\bibitem[SNCH04]{smith2004multiscale}
N.P. Smith, D.P. Nickerson, E.J. Crampin, and P.J. Hunter.
\newblock Multiscale computational modelling of the heart.
\newblock {\em Acta Numerica}, 13:371--431, 2004.

\bibitem[St{\"u}01]{stuben2001review}
K.~St{\"u}ben.
\newblock A review of algebraic multigrid.
\newblock {\em Numerical Analysis: Historical Developments in the 20th
  Century}, pages 331--359, 2001.

\bibitem[TW04]{toselli2004domain}
A.~Toselli and O.~Widlund.
\newblock {\em Domain decomposition methods-algorithms and theory}, volume~34.
\newblock Springer Science \& Business Media, 2004.

\bibitem[ULM02]{Usyk2002}
T.P. Usyk, I.J. LeGrice, and A.D. McCulloch.
\newblock Computational model of three-dimensional cardiac electromechanics.
\newblock {\em Computing and Visualization in Science}, 4(4):249--257, Jul
  2002.

\bibitem[WN99]{wright1999numerical}
S.~Wright and J.~Nocedal.
\newblock Numerical optimization.
\newblock {\em Springer Science}, 35(67-68):7, 1999.

\bibitem[XZ17]{xu2017algebraic}
J.~Xu and L.~Zikatanov.
\newblock Algebraic multigrid methods.
\newblock {\em Acta Numerica}, 26:591--721, 2017.

\bibitem[Zam16]{zampini2016pcbddc}
S.~Zampini.
\newblock {PCBDDC}: a class of robust dual-primal methods in {PETSc}.
\newblock {\em SIAM Journal on Scientific Computing}, 38(5):S282--S306, 2016.

\end{thebibliography}
